%% file: lagrangian_14_giugno_2010_.tex
\documentclass[reqno,11pt]{amsart}

\usepackage{amsmath,amsfonts,amssymb,amsthm,epsfig,amscd}
\usepackage[]{fontenc}
\usepackage{xy}
\usepackage{enumerate}

 \allowdisplaybreaks \numberwithin{equation}{section}

\def\stepempty#1{(#1)}

\xyoption{all} \numberwithin{equation}{section}

\newtheorem{theorem}{Theorem}[section]
\newtheorem{proposition}[theorem]{Proposition}
\newtheorem{conjecture}[theorem]{Conjecture}
\newtheorem{lemma}[theorem]{Lemma}

\theoremstyle{definition}
\newtheorem{remark}[theorem]{Remark}
\newtheorem{definition}[theorem]{Definition}
\newtheorem{example}[theorem]{Example}

\newtheorem{notation}[theorem]{Notation}

\newcommand{\C}{\mathbb{C}}
\newcommand{\Q}{\mathbb{Q}}
\newcommand{\Z}{\mathbb{Z}}

\newcommand{\pr}{\mathbb P}
\newcommand{\sym}{\mbox{\upshape{Sym}}}

\newcommand{\om}{\omega}

\newcommand{\oo}{\mathcal{O}}
\newcommand{\pic}{\mbox{\upshape{Pic}}}

\newcommand{\ba}{{B_{\alpha_1}}}
\newcommand{\rg}{{R_{\overline \gamma}}}
\newcommand{\prym}{{P(\pi)}}
\newcommand{\Alb}{\mbox{Alb}}
\newcommand{\alb}{{a}}

\newcommand{\bC}{{\mathbb C}}

\newcommand{\lpr}{LG}

\begin{document}

\title[Galois closure and lagrangian varieties]{Galois closure and lagrangian varieties}
\date{December 11, 2009}
\author{F. Bastianelli}\address{Dipartimento di Matematica, via Ferrata 1, 27100 Pavia}\email{francesco.bastianelli@unipv.it}
\author{G. P. Pirola}\address{Dipartimento di Matematica, via Ferrata 1, 27100 Pavia}\email{gianpietro.pirola@unipv.it}
\author{L. Stoppino}\address{Dipartimento di Fisica e Matematica, via Valleggio 11, 22100 Como}\email{lidia.stoppino@uninsubria.it}

\thanks{This work has been partially supported by 1) PRIN 2007 \emph{``Spazi di moduli e teorie di Lie''}; 2) Indam (GNSAGA); 3) FAR 2008 (PV) \emph{``Variet\`a algebriche, calcolo algebrico, grafi orientati e topologici''.}}

\maketitle

\begin{abstract}
Let $X$ be a complex projective variety and consider the morphism
$$
\psi_k\colon \bigwedge^k H^0(X, \Omega^1_X)\longrightarrow H^0(X, \Omega^k_X).
$$
We use Galois closures of finite rational maps to introduce a new method for producing varieties such that $\psi_k$ has non-trivial kernel.
We then apply our result to the two-dimensional case and we construct a new family of surfaces which are Lagrangian in their Albanese variety. Moreover, we analyze these surfaces computing their Chern invariants, and proving that they are not fibred over curves of genus $g\geq 2$. The topological index of these surfaces is negative and this provides a counterexample to a conjecture on Lagrangian surfaces formulated in \cite{BNP}.
\end{abstract}

\section*{Introduction}

Let $X$ be a smooth complex projective variety of dimension $n$ and let ${H^\ast (X, \mathbb C)=\bigoplus_{k\geq 0} H^k(X,\mathbb C)}$ be  its complex cohomology.
Since the fundamental paper \cite{DGMS}, it has been known that the algebra structure, induced by the cup product, determines the rational homotopy of $X.$
We consider the exterior algebra $ \bigwedge^\ast H^1(X,\mathbb C)=\bigoplus_{k\geq 0}\bigwedge^kH^1(X,\mathbb C)$
and the natural homomorphism
$$
\rho\colon \bigwedge^*H^1(X,\mathbb C)\equiv H^*(\Alb(X), \mathbb C)\longrightarrow H^*(X, \mathbb C),
$$
which can be seen as the algebra homomorphism induced from the Albanese map $\alb \colon X\longrightarrow \Alb(X).$
It follows that  $H^\ast (X,\mathbb C)$ is a graded $\bigwedge^\ast H^1(X,\mathbb C)$-module.
This structure also provides important information on the topology of $X$.
A basic result in this direction (see  \cite{Morgan1}, \cite{Hein} , \cite{Campana} and \cite{wonderful}) is that the nilpotent completion of the fundamental group $\pi_1(X)$ is determined by $H^1(X,\Q)$ and by the homomorphism
$$
\rho_2\colon \bigwedge^2 H^1(X,\C) \longrightarrow H^2(X,\C).
$$

In particular  if $\ker\rho_2$ is not trivial, then  $\pi_1(X)$ is not abelian.
Examples of varieties with non-trivial nilpotent towers were found by Sommese-Van de Ven \cite{SVdV} and by Campana \cite{Campana}.
In these examples the non-trivial elements in $\ker \rho_2$ all come from the $(1,1)$-part of $\rho_2$ in the Hodge decomposition
$\rho_2^{1,1}\colon H^{1,0}(X)\otimes H^{0,1}(X)\longrightarrow H^{1,1}(X)$.

\smallskip
In this paper we focus instead on the holomorphic part of $\rho$.
Consider the algebra of holomorphic forms on $X$
$$H^{\ast , 0}(X)=\bigoplus_{k\geq 0}H^{k,0}(X)= \bigoplus_{k\geq 0}H^0(X, \Omega^k_X),$$
and on the Albanese variety $\Alb(X) $
$$
H^{\ast , 0}(\mbox{Alb}(X))=\bigwedge^\ast H^{1,0}(X)= \bigoplus_{k\geq 0}\bigwedge^k H^{1,0}(X) =\bigoplus_{k\geq 0} \bigwedge^k H^0(X, \Omega^1_X).
$$
The holomorphic part of $\rho$
\begin{equation}\label{G}
\psi:=\rho^{\ast,0}\colon \bigwedge^\ast H^{1,0}(X)\longrightarrow H^{\ast , 0}(X)
\end{equation}
provides $H^{\ast , 0}(X)$ of a structure of $\bigwedge^\ast H^{1,0}(X)$-module.

\smallskip
The importance of these structures  has been emphasized recently by  Lazarsfeld and Popa in \cite{L-P}, where generic vanishing is applied to the setting of the BGG correspondence.

A classical result relating  the kernel of $$ \psi_2\colon \bigwedge^2 H^{1,0}(X)\longrightarrow H^{2,0}(X)$$ to the topology of $X$ is the Castelnuovo-de Franchis  Theorem (cf. \cite[Proposition X.9]{Beauville}).
It states that  the  existence of  $w_1, w_2 \in H^0(X, \Omega^1_X)$   such that $0\not =w_1\wedge w_2\in\ker\psi_2$  is equivalent to the existence of a surjective morphism from $X$  to a curve of genus $\geq 2$.
We note that in this case the fundamental group of $X$  surjects onto the fundamental group of the base curve; hence it has an infinite nilpotent tower, in line with the above mentioned results.
The Castelnuovo-de Franchis Theorem has been generalized to the case of decomposable elements in  $\ker \psi_k$ for arbitrary $k$ by Catanese \cite{Catanese}.

Apart from the decomposable case, the most important setting where non-trivial elements appear in $\ker \psi_2$ is when $X$ has a Lagrangian structure.
We will say that an $n$-dimensional variety  $X$  is {\em Lagrangian} if there exists a map of degree one $\alpha\colon X\longrightarrow \alpha(X)\subset A$, where $A$ is an abelian variety of dimension $2n$, and a $(2,0)$-form $w$ of rank $2n$ on $A$ such that $\alpha^*w=0$. Clearly, a Lagrangian subvariety of an abelian variety is Lagrangian.
One expects that the existence of a Lagrangian form gives strong restrictions on the geometry and the topology of a variety; however, it turned out that even to construct examples of such varieties is a difficult task.

In  \cite{BT} Bogomolov and Tschinkel provide examples of surfaces which are Lagrangian in their Albanese variety, using dominant maps between K3 surfaces. In some cases they prove that these surfaces are not fibred over curves of genus $g\geq 2$.

In  \cite{BNP} a new topological consequence of the existence of a Lagrangian form is proved; namely, if $X$ is a Lagrangian surface, under some assumption  on the branch locus of the Albanese map, the topological index $\tau(X)=\frac{1}{3}(K_X^2-2c_2(X))$ is showed to be non-negative. Moreover, the authors give another example of a surface with a Lagrangian $(2,0)$-form not coming from a fibration (\cite[Example 6.6]{BNP}) based on a previous article of the third author \cite{Pirola}.

\bigskip

In this paper we investigate the map $\psi$ for a variety $X$ constructed as Galois closure of a generically finite morphism $\gamma \colon Z\longrightarrow Y$ (or else a finite degree dominant rational map). We study the vector spaces $H^{p,0}(X)$ and $ \bigwedge^p H^{1,0}(X)$ as representations of the Galois group of $\gamma$. Recall that such a group is isomorphic to the monodromy group of $\gamma$. We suppose that $Z$ is irregular - i.e. $h^{1,0}(Z)>0$ - and that ${h^{1,0}(Y)=h^{p,0}(Y)=0}$ for some $p$. By using representation theory we then detect a subspace of large dimension contained in $\ker \psi_p$.

Although the Galois closures have been thoroughly studied by several authors (see for instance the important papers \cite{Miyaoka} and \cite{M-T}), this very simple and general application of the theory of Galois covers is - to our knowledge - new.

We can summarize our result as follows (see Theorem \ref{theorem KERNEL} for a stronger statement).
\begin{theorem}\label{teo1}
Let $\gamma \colon Z\dasharrow Y$ be a  finite degree $d$ map between smooth varieties of dimension $n$ and let $g\colon X\dasharrow Y$ be the Galois closure of $\gamma$. Suppose that
\begin{itemize}
  \item[(i)] the monodromy group $M(\gamma)$ is isomorphic to the full symmetric group $S_d$,
  \item[(ii)] the irregularity of $Z$ is $q=h^{1,0}(Z)\geq 2$,
  \item[(iii)] $h^{1,0}(Y)=0$ and there exists an integer $2\leq p\leq q$ such that $h^{p,0}(Y)=0$.
\end{itemize}
Then the kernel of  $\psi_p\colon \bigwedge^pH^{1,0}(X)\longrightarrow H^{p,0}(X)$ has dimension grater or equal than ${q \choose p}$.
\end{theorem}
The main point is to describe the representation of the group $S_d$ given by $H^{1,0}(X)$ to find a trivial sub-representation in the decomposition of $\bigwedge^p H^{1,0}(X)$. Therefore this subspace has to lie into $\ker \psi_p$, because of assumption (iii).

We note further that - under the assumption of the theorem - we are able to exhibit explicitly ${q \choose p}$ independent elements of the kernel of $\psi_p$ (cf. Proposition \ref{proposition p-FORMS}).

In \cite[Section 3]{Campana} Campana underlines the importance of producing examples of varieties such that the holomorphic $(2,0)$-part of $\rho$ has non-trivial kernel. Observe that the above theorem allows to construct many varieties with such a property. In particular, the fundamental groups of these varieties have non-trivial nilpotent tower.

\medskip
In the second part of the paper we construct and study in detail a family of Lagrangian surfaces of general type, which we call $\lpr$ \emph{surfaces} (cf. Definition \ref{definition LPR}). These surfaces are Galois covers of a triple covering of an abelian surface on a rational one, so the existence of a non-trivial element in $\ker \psi_2$ is guaranteed by Theorem \ref{theorem MAIN}. On the other hand, the study of the geometry of $\lpr$ surfaces requires a lot of effort and the use of several different techniques. Besides classical tools (such as linear systems on curves, fibrations and double coverings on surfaces), we exploit the theory of abelian surfaces and their moduli, and the monodromy of the torsion points.

The starting point of our construction is an abelian surface $S$ provided of a line bundle $\mathcal L$ inducing a polarization of type $(1,2)$. Using the fine analysis of $(1,2)$-polarized abelian surfaces contained in \cite{Barth}, we investigate the geometry of the linear pencil $|\mathcal L|$ and of the induced fibration on $\pr^1$ relating the behavior of the special fibers to certain subsets of the torsion points of $S$.
We then define a degree three rational covering $\gamma\colon S\dashrightarrow \mathbb F_3$ of the Hirzebruch surface $\mathbb F_3$. We would like to note that the induced fibration has been studied also by Xiao \cite{Xiao} and by Chen and Tan in \cite{CT}, whereas the map $\gamma$ has been presented by Tokunaga and Yoshihara \cite{TY} in order to provide examples of abelian surfaces with minimal degree of irrationality.

We give an explicit geometric description of the Galois closure of $\gamma$, and we define the $\lpr$ surface $X$ to be its minimal desingularization. Even though this construction could be carried on a birational setting, in order to control the geometry of $X$ and to compute its invariants, we need to resolve explicitly the indeterminacy locus of $\gamma$: it is then crucial the knowledge of the geometry of $|\mathcal L|$.

We investigate the geometry of $\lpr$ surfaces.
The main results we obtain can be summarized as follows.
\begin{theorem}\label{theorem MAIN}
Let $S$ be an abelian surface with a line bundle $\mathcal L$ of type $(1,2)$ such that all the curves in the pencil $|\mathcal L|$ are irreducible. Then there exists a rational triple covering $\gamma \colon S\dashrightarrow \mathbb{F}_3$ of the Hirzebruch surface $\mathbb{F}_3$. The minimal desingularization $X$ of the Galois closure of $\gamma$ is a surface of general type with invariants
$$
K_{X}^2= 198,\quad c_2(X)=102,\quad \chi(\oo_{X})=25,\quad q=4 ,\quad p_g=28\quad and\quad  \tau(X)=-2.
$$
Moreover, $X$ is Lagrangian with Albanese variety  $\Alb (X)=S\times S$, and it does not have any fibration over curves of genus $\geq 2$.
\end{theorem}
We prove further that for a general choice of $S$ in the moduli space $\mathcal W(1,2)$ of $(1,2)$-polarized abelian surfaces, the Galois closure of the covering $\gamma$ is smooth and minimal, hence it is the $\lpr$ surface itself. A crucial point here is the study of the monodromy of the points of order $3$ in $S.$

In view of Castelnuovo-de Franchis Theorem the $\lpr$ surfaces are particularly interesting because the Lagrangian structure does not come from a fibration over a curve of genus greater or equal to $2.$ Our  proof of this fact is indirect and uses heavily the knowledge of special curves on the $\lpr$ surface and the theory of complex abelian surfaces.

It is worth noticing that the $\lpr$ surfaces are so far the only examples of  Lagrangian surfaces non-fibred over curves of genus greater than $1$ whose invariants are explicitly computed.

\medskip

We see from Theorem \ref{theorem MAIN} that the $\lpr$ surfaces have negative topological index $\tau(X)=-2$.
This shows that the positivity of the index is not a property of Lagrangian surface, thus disproving a conjecture stated in \cite{BNP}.
The assumptions on the branch locus of the Albanese map made in that paper, and in particular the connectedness,  are then proved to be necessary.
In Section \ref{topological} we make a detailed discussion on this topic.

\medskip

From the non-triviality of $\ker \rho_2$ it follows straightforwardly that  the nilpotent tower of the $\lpr$ surfaces is not trivial up to the second step. The examples of Campana and Sommese-Van der Ven have precisely two steps.
In \cite{ATV} (the last paper of a series) Amram, Teicher and Vishne study the Galois closures of  the generic projections from the product of two elliptic curves to $\mathbb P^2$. Among other results, they prove that the nilpotent tower of the fundamental groups of these surfaces are non-trivial up to step $3$.
It would be interesting to compute the class of nilpotency of the $\lpr$ surfaces..  
We determine the dimension of $\ker\rho_2$ (see Proposition \ref{kerrho2}) which is $7$. This is a first step towards the computation of the nilpotent tower of $\pi_1(X)$, which is one of our next purposes.


\subsection*{Plan of the paper} In the first section we prove the general result described in Theorem \ref{teo1}. The rest of the paper is devoted to the construction and study of the $\lpr$ surfaces. In the second section we recall some results on $(1,2)$-polarized abelian surfaces, and we establish a correspondence between bielliptic curves and these surfaces (via generalized Prym varieties). In Section 3 we consider the pencil defined by the polarization on one of these abelian surfaces $S$, and we relate the geometry of special members to some subsets of torsion points on $S$. We define in Section 4 a triple covering from a blow up $\overline S$ of $S$ to the Hirzebruch surface $\mathbb F_3$, whose Galois closure $\pi_1\colon X\longrightarrow \overline S$ is  constructed in Section 5. Moreover, we prove that the branch locus $B_{\alpha_1}$ is reduced with simple double singularities and we  compute its numerical class. In Section 6 we then define $\lpr$ surfaces and we compute their Chern invariants. We determine their Albanese varieties and prove that they are Lagrangian. The seventh section is devoted to prove that the $\lpr$ surfaces do not admit a fibration over curves of genus $g\geq 2$, and we study $\ker\rho_2$. In Section 8, via the monodromy of points of order $3$, we prove that the Galois closures associated to general abelian surfaces are smooth. This result is crucial for the study of the branch locus of the Albanese map carried on in Section 9. Here we relate our results to the ones of \cite{BNP}.

\subsection*{Notation}

We work over the field $\mathbb C$ of complex numbers. By \emph{variety} we shall mean a complex complete irreducible algebraic variety, unless otherwise stated.

Given a smooth variety $X$, we denote by $\Omega^1_X$ the cotangent bundle of $X$, by $\Omega^{p}_X=\bigwedge^p \Omega^1_X$ the sheaf of holomorphic $p$-forms and by $h^{p,q}(X)=\dim H^{q}(X,\Omega^{p})$ the Hodge numbers, where $p$ and $q$ are non-negative integers; in particular $q(X)=h^{1,0}(X)$ is the irregularity of $X$. We recall that $h^{p,0}(X)$ are birational invariants. We denote by $\omega_X$ the canonical sheaf of $X$, which is the line bundle $\omega_X=\bigwedge^n \Omega^1_X$ of holomorphic $n$-forms, with $n=\dim X$. Moreover, we denote by $K_X$ any divisor such that $\mathcal O_X(K_X)\cong \omega_X$.

The Albanese variety of $X$ is the abelian variety $\Alb(X)=H^{0}(X, \Omega^1_X)^\vee/H_1(X,\Z)$. Given a point $\overline x\in X$ we can define the Albanese map $\alb_{\overline x}\colon X\longrightarrow \Alb(X)$ by $x\mapsto \int_{\overline x}^x-$. Changing the base point the Albanese map changes by a translation of $\Alb(X)$. We will therefore often omit the base point.

If $D_1$ and $D_2$ are two Weil divisors on $X$, the notation $D_1\sim_X D_2$ means that $D_1$ and $D_2$ are linearly equivalent, i.e. that there is an isomorphism of line bundles $\mathcal O_X(D_1)\cong \mathcal O_X(D_2)$. By $D_1\equiv_X D_2$ we mean that the two divisors are numerically equivalent, that is $D_1\cdot C=D_2\cdot C$ for any curve $C\subset X$.

A {\em fibration} $f\colon X\longrightarrow Y$ over a smooth projective variety $Y$ is a flat proper surjective morphism with connected fibers. We denote by $\omega_f$ the relative canonical sheaf $\omega_f=\omega_X\otimes f^*\omega_Y^{-1}$.

We say that a property holds for a \emph{general} point $x\in X$ if it holds on an open non-empty subset of $X$. Moreover, we say that $x\in X$ is a \emph{very general} - or \emph{generic} - point if there exists a countable collection of proper subvarieties of $X$ such that $x$ is not contained in the union of those subvarieties.

\subsection*{Acknowledgments} We would like to thank Jaume Amor\'os and Fr\'ed\'eric Campana for helpful suggestions.

\medskip
\section{Galois theory and irregularity}\label{section GALOIS-IRREGULARITY}

Let $Z$ and $Y$ be two smooth varieties of dimension $n$ and let $\gamma\colon Z \dasharrow Y$ be a degree $d>1$ rational map. Let $K(Y)$ and $K(Z)$ denote the fields of rational functions. Then $\gamma$ induces a field extension $K(Y)\subseteq K(Z)$ of degree $d$. Let $L$ be the Galois closure of this field extension and let $W$ be the normalization of $Z$ in $L$ (see for instance \cite{Liedtke}).
\begin{definition}\label{definition GALOIS CLOSURE}
With the above notation, the normal variety $W$ provided of the induced map $g\colon W\dasharrow Y$ is the \emph{Galois closure} of the map $\gamma$.
\end{definition}
Let $X$ be a desingularization of $W$ and, for $p\geq 2$, let us consider the homomorphism
\begin{equation}\label{equation CUP PRODUCT MAP}
\psi_p \colon \bigwedge^p H^{1,0}(X)\longrightarrow H^{p,0}(X).
\end{equation}
In the following we  study of the kernel of this map and prove that - under some additional hypothesis - it is non-trivial. Note that $\psi_p$ is independent from the desingularization chosen.  We can therefore suppose - by blowing up its indeterminacy locus - that $\gamma$ is a generically finite  morphism of degree $d$.\\

Let $M=M(\gamma)$ be the monodromy group of $\gamma$ \cite{Harris} and let $m$ be its order. Recall that it is isomorphic to the Galois group of the splitting normal extension $L/K(Y)$ induced by $\gamma$. Moreover, the Galois group of the extension $K(W)/K(Z)$ is a subgroup $M_1\subset M$.

It is possible to describe geometrically the Galois closure $W$ of $\gamma$ as follows (cf. \cite[Expos\'e V.4.g]{SGA}, \cite{M-T} and \cite[Proposition 6.13]{Liedtke}). Let $\mathcal W\subset Y$ be a suitable non-empty Zariski open subset of $Y$ and let $\mathcal{U}:=\gamma^{-1}(\mathcal{W})\subset Z$ so that the restriction $\gamma_{| \mathcal U}\colon \mathcal U \longrightarrow \mathcal W$ is an \'etale morphism of degree $d$. Let $Z^d=Z\times\ldots\times Z$ be the $d$-fold ordinary product of $Z$ and let us consider the subset $W^\circ\subset Z^d$ given by
$$
W^\circ :=\left\{\left.(z_1,\ldots,z_d) \in Z^d\right| \{ z_1,\ldots,z_d \}=\gamma^{-1}(y) \textrm{ for some } y\in \mathcal W \right\}.
$$
Clearly $W^\circ$ is a smooth - possibly disconnected - variety and the full symmetric group $S_d$ acts on $W^{\circ}$. By lifting arcs it is immediate to see that $W^\circ$ has exactly $\frac{d!}{m}$ irreducible components, that are all isomorphic. We recall that $d$ divides $m$ and that $m=d$ if and only if $\gamma$ is a Galois covering. Furthermore, $W^\circ$ is irreducible if and only if $M\cong S_d$.

Let us fix an irreducible component of $W^\circ$ and let $W$ be the normalization of the Zariski closure in $Z^d$ of such a component. Let $\alpha_i\colon W\longrightarrow Z$ be the composition of the normalization morphism with the $i$-th projection map from $Z^d$ to $Z$, with $i=1,\ldots,d$. Thus we have the following diagram
$$
\xymatrix{ W\ar[dr]^{g}\ar[d]_{\alpha_1} & \\ Z\ar[r]^{\gamma} & Y\\}
$$
\smallskip
and the normal variety $W$ - provided with the morphism  $g=\gamma\circ \alpha_1\colon X\longrightarrow  Z$ - turns out to be the Galois closure of the covering $\gamma\colon Z\longrightarrow Y$. In particular, the monodromy group of the map $\alpha_1$ is $M(\alpha_1)=M_1\subset M:=M(\gamma)$.
\begin{remark}\label{remark Z^(d-1)}
We note that the variety $W^\circ$ can also be defined as the subvariety of $Z^{d-1}$
\begin{displaymath}
W^\circ :=\left\{\left.(z_1,\ldots,z_{d-1}) \in Z^{d-1}\right| \{ z_1,\ldots,z_{d-1},z_d \}=\gamma^{-1}(y) \textrm{ for some } y\in \mathcal W \right\}.
\end{displaymath}
The morphisms $\alpha_1,\ldots,\alpha_{d-1}\colon W\longrightarrow Z$ are defined by composing the normalization map $W\longrightarrow W^\circ$ with the projection maps, whereas $\alpha_d\colon W\longrightarrow Z$ is obtained  sending $(z_1,\ldots,z_{d-1})\in Z^{d-1}$ to the point $z_d\in Z$.
\end{remark}

Now, let $\mu\colon X\longrightarrow W$ be a desingularization of $W$ and, for $i=1\ldots,d$, let us consider the maps $\pi_i:=\alpha_i\circ\mu\colon X\longrightarrow Z$. We have the following diagrams
\begin{equation}\label{diagram XZY}
\xymatrix{ X\ar[dr]\ar[d]_{\pi_i} & \\ Z\ar[r]^{\gamma} & Y\\}
\end{equation}
and the monodromy groups are $M(\pi_1)=M(\alpha_1)=M_1$ and $M(\gamma\circ\pi_1)=M(\gamma)=M$. We remark that $M$ acts birationally on $X$ and hence the vector space $H^{k,0}(X)$ carries - in a natural way - a $M$-representation structure for any $k\geq 0$. Moreover, we have the following relations on the spaces of holomorphic differential forms:
\begin{equation}\label{equation RELATIONS}
H^{k,0}(X)^{M}=  (\gamma\circ\pi_1)^{\ast}(H^{k,0}(Y))\quad \mbox{ and } \quad H^{k,0}(X)^{M_1}= \pi_1^{\ast}(H^{k,0}(Z)).
\end{equation}
We note that the cup product map
$$
H^{k,0}(X)\otimes H^{h,0}(X)\longrightarrow H^{k+h,0}(X)
$$
is $M$-equivariant for any $k,h\geq 0$. Furthermore, the map $\psi_k\colon \bigwedge^k H^{1,0}(X)\longrightarrow H^{k,0}(X)$ introduced in (\ref{equation CUP PRODUCT MAP}) is $M$-equivariant as well.\\

For the sake of simplicity we assume hereafter that the monodromy group of $\gamma$ is
the full symmetric group  $S_d$.
In particular, it follows that $X$ is the normalization of the Zariski closure of $W^\circ$ and the monodromy group of $\pi_1$ is $M_1\cong S_{d-1}$.

Recall that the \emph{standard representation} of $S_d$ is the irreducible representation defined as the $(d-1)$-dimensional $\mathbb{C}$-vector space $\Gamma:=\left\{(a_1,\ldots,a_d)\in \mathbb{C}^d|a_1+\ldots+a_d=0\right\}$, where we fixed the standard basis $\left\{e_1,\ldots,e_d\right\}$ of $\mathbb{C}^d$ and $\sigma e_i=e_{\sigma(i)}$ for any $\sigma\in S_d$, $i=1,\ldots,d$. Moreover, let $U$ denote the \emph{trivial representation} of $S_d$. Then we have the following.
\begin{theorem}\label{theorem qCOPIE}
Suppose that $h^{1,0}(Y)=0$. Then the vector space $H^{1,0}(X)$ contains $q:=q(Z)$ copies of the standard
representation $\Gamma$ of $S_d$. In particular,
\begin{equation}\label{equation h1,0}
h^{1,0}(X)\geq q(d-1).
\end{equation}
\begin{proof}
Let us consider the canonical identification  $H^{1,0}(Z^d)\cong H^{1,0}(Z)^d$. Let $\omega \in H^{1,0}(Z)$ be a non-zero holomorphic $1$-form and let us set $\omega^i:=\pi_i^{\ast}(\omega)\in H^{1,0}(X)$ for any $i=1,\ldots,d$. We note that the $\omega^i$'s are distinct non-zero forms on $X$. Indeed, if it were $\omega^i=\omega^j$ for some $i\neq j$, then $\omega^j$ would be invariant with respect to both ${M(\pi_i)\cong S_{d-1}}$ and ${M(\pi_j)\cong S_{d-1}}$. Hence the non-zero form $\omega^j$ would be invariant under the action of the whole $S_d$, but this is impossible because $Y=X/S_d$ does not posses holomorphic $1$-forms.

Moreover, as the $1$-form $\sum_i\omega^i$ is invariant under the action of $S_d$ and $h^{1,0}(Y)=0$, we have that ${\omega^1+\ldots+\omega^d=0}$. Thus the vector space $\langle \omega^1, \ldots , \omega^d\rangle$ - carrying a $S_d$-representation structure - consists of a copy $\Gamma_{\omega}$ of the standard representation $\Gamma$ of $S_d$. Let
$$
\eta \colon H^{1,0}(Z)\otimes \Gamma \longrightarrow H^{1,0}(X)
$$
be the homomorphism sending an element $\omega \otimes a$, with $\omega\in H^{1,0}(Z)$ and $a=(a_1,\ldots,a_d)\in \Gamma$, to the corresponding element $\omega_a=\sum_i a_i \omega^i$.

To conclude the proof it remains to show that the homomorphism $\eta$ is injective.
Note that $\eta$ is $M$-equivariant, where the action of $M=S_d$ on $H^{1,0}(Z)$ is the trivial one.
Thus the kernel of $\eta$ is a sub-representation of $H^{1,0}(Z)\otimes \Gamma$ and hence $\ker(\eta)=\Sigma \otimes  \Gamma$, where $\Sigma \subset H^{1,0}(Z)$ is a trivial sub-representation. Furthermore, we have that
$$
(H^{1,0}(Z)\otimes \Gamma)^{M_1}=(H^{1,0}(Z)\otimes \Gamma)^{S_{d-1}}\cong  H^{1,0}(Z)\otimes \Gamma^{S_{d-1}}\cong H^{1,0}(Z)
$$
and then $\ker(\eta )^{M_1}=\Sigma$. Since $\eta^{M_1}\colon (H^{1,0}(Z)\otimes \Gamma)^{M_1}\longrightarrow H^{1,0}(X)^{M_1}$ coincides with the injective homomorphism $\pi_1^{\ast}\colon H^{1,0}(Z)\longrightarrow H^{1,0}(X)$, we deduce that $\Sigma=\ker(\eta )^{M_1}=\ker(\eta ^{M_1})=\{0\}$. We conclude that $\ker(\eta)=\{0\}$.
\end{proof}
\end{theorem}

\begin{remark}
We note that  inequality (\ref{equation h1,0}) can be strict. For instance, let us consider the case of curves. As is  well known, a smooth curve $C$ of genus $g$ admits a degree $d$ covering $\gamma\colon C\longrightarrow \mathbb{P}^1$ with monodromy group $M(\gamma)=S_d$ for $d\geq\frac{g+2}{2}$.
Let $X$ be the Galois closure of such a covering and assume that $d\geq 6$. Then $X$ admits a covering $X\longrightarrow C$ of degree $(d-1)!$ and by Riemann-Hurwitz formula we have
$$
h^{1,0}(X)\geq(d-1)!\,(g-1)+1>(d-1)!\,(q-1)\geq (d-1)q.
$$
\end{remark}

\bigskip
Thanks to Theorem \ref{theorem qCOPIE} we have that $H^{1,0}(X)$ contains $\Gamma^{\oplus q}$ as a sub-representation. Hence if $q\geq 2$ and $p$ is an integer such that $2\leq p\leq q$, the $S_d$-representation $\bigwedge^p H^{1,0}(X)$ admits a decomposition of the form
$$
\bigwedge^p H^{1,0}(X)\cong \left(\otimes^{p}\,\Gamma\right)^{\oplus r} \oplus K,
$$
where $r:={q \choose p}$ and $K$ is some sub-representation. It follows that
\begin{equation}\label{equation SYMp INTO WEDGEp}
\left(\sym^p\,\Gamma\right)^{\oplus r}\subset \bigwedge^p H^{1,0}(X)
\end{equation}
is a sub-representation. In order to study the kernel of the map $\psi_p$ defined in (\ref{equation CUP PRODUCT MAP}), we shall discuss the existence of trivial representations of $S_d$ into $\bigwedge^p H^{1,0}(X)$. Then we are going to focus on the trivial representations contained in $\sym^p\,\Gamma$.

Letting  $V=U\oplus \Gamma$, the graded algebra $\oplus_k \sym^k \,V$ can be naturally identified with the algebra of polynomials over $\mathbb{C}$ in the variables $x_1,\ldots,x_d$, that is
$$
S(x_1,\ldots,x_d)\cong \oplus_k \sym^k \,V.
$$
The ring of the $S_d$-invariant polynomials of $S(x_1,\ldots,x_d)$ is generated by the elementary symmetric functions
\begin{equation}\label{equation EL SYMM FUN}
\xi_h(x_1,\ldots,x_d):=\sum_{1\leq i_1<i_2<...<i_h\leq d} x_{i_1}x_{i_2}...x_{i_h},
\end{equation}
with $1\leq h\leq d$. Therefore
$$
S(x_1,\ldots,x_d)=I(\xi_1)\oplus \left(\oplus_k \sym^k\,\Gamma\right),
$$
where $I(\xi_1)$ is the principal ideal generated by $\xi_1$. Then the sub-algebra of the $S_d$-invariant elements of $\oplus_k \sym^k\,\Gamma$ can be identified with the $\mathbb{C}$-algebra generated by the functions $\xi_2,\ldots,\xi_d$, i.e. we have the isomorphism of graded algebras given by
\begin{equation}\label{equation IDENTIFICATION}
\left(\oplus_k \sym^k\,\Gamma\right)^{S_d}\cong \mathbb{C}[\xi_2,\ldots,\xi_d].
\end{equation}
Summing up, we proved the following.
\begin{lemma}\label{lemma INVARIANT ELEMENTS OF SYMkGAMMA}
Let $\Gamma$ be the standard representation of $S_d$ and let $\xi_1,\ldots,\xi_d$ be the elementary symmetric polynomials defined in (\ref{equation EL SYMM FUN}). Then the vector space
\begin{equation}\label{equation Ak}
A_k:=\bC[\xi_2,\ldots,\xi_d]_k
\end{equation}
is isomorphic to the space of the $S_d$-invariant elements of $\sym^k\,\Gamma$ for any $k\geq 0$.
\end{lemma}

In particular, every $A_k\cong \left(\sym^k\,\Gamma\right)^{S_d}$ is a sum of copies of the trivial representation of $S_d$. Moreover, by (\ref{equation SYMp INTO WEDGEp}) we have
\begin{equation}\label{equation Apr NEL WEDGEp}
 (A_p)^{\oplus r}\subset \left(\bigwedge^p H^{1,0}(X)\right)^{S_d},
\end{equation}
where $2\leq p\leq q$ and $r={q \choose p}$, i.e. $\bigwedge^p H^{1,0}(X)$ contains $r$ copies of $A_p$ as trivial sub-representations. We are now ready to state the main result of this section.
\begin{theorem}\label{theorem KERNEL}
Let $\gamma \colon Z\dasharrow Y$ be a finite map of degree $d$ between smooth varieties of dimension $n$ and let $g\colon X\dasharrow Y$ be the Galois closure of $\gamma$. Suppose that
\begin{itemize}
  \item[(i)] the monodromy group $M(\gamma)$ is isomorphic to the full symmetric group $S_d$,
  \item[(ii)] $q=q(Z)\geq 2$,
  \item[(iii)] $h^{1,0}(Y)=0$ and there exists an integer $2\leq p\leq q$ such that $h^{p,0}(Y)=0$.
\end{itemize}
Then the cup-product map
$$
\psi_p\colon \bigwedge^p H^{1,0}(X)\longrightarrow H^{p,0}(X)
$$
has non-trivial kernel and moreover
\begin{equation}\label{equation Ap nel ker}
(A_p)^{\oplus r}\subset \ker \psi_p\,,
\end{equation}
where $r={q \choose p}$ and $A_p$ is the vector space defined in (\ref{equation Ak}). In particular,
$$
\dim\ker \psi_p \geq {q \choose p}.
$$
\begin{proof}
As $H^{p,0}(Y)=\{0\}$, by (\ref{equation RELATIONS}) we have that $\{0\}=g^{\ast}(H^{p,0}(Y))=H^{p,0}(X)^{S_d}$. Clearly, $\psi_p\left(\bigwedge^p H^{1,0}(X)\right)^{S_d}\subset H^{p,0}(X)^{S_d}$ because $\psi_p$ is $S_d$-equivariant. Hence $(A_p)^{\oplus r}\subset \ker \psi_p$ by (\ref{equation Apr NEL WEDGEp}). In particular, we deduce that $\dim\ker \psi_p\geq {q \choose p}\dim A_p\geq {q \choose p}$.
\end{proof}
\end{theorem}

\begin{remark}
In order to construct non-trivial elements in the kernel of $\psi_p$, the assumptions (i) and (iii) can be weakened. On one hand, it is not necessary to suppose that the monodromy group is the full symmetric group and it suffices to assume that $\gamma$ is not a Galois covering (for instance one could consider the case of the alternating group $A_d$ with $d\geq 4$).
On the other hand, if we drop the assumption $q(Y)=0$, we have that $\ker \psi_p$ contains ${a \choose p}$ copies of $A_p$ with $a:=q(Z)-q(Y)$.
\end{remark}
\begin{remark}
It is worth noticing that Theorem \ref{theorem KERNEL} is more significant when $p\leq \dim \alb(Z)$. Indeed, for $p$ strictly grater than the Albanese dimension of $Z$ it is easy to construct subspaces of $\ker \psi_p$. Consider the Albanese maps of $X$ and $Z$, and let $h\colon \Alb (Z)\longrightarrow \Alb (X)$ be the morphism induced from $g\colon X\longrightarrow Z$. We have the commutative diagram induced in cohomology
$$
\xymatrix{ \quad\bigwedge^pH^{1,0}(X)=  \!\!\!\!\!\!\!\!\!\!\!\!\!\!\!\!&H^{p,0}(\Alb(X))\ar[r]^{\psi^X_p}& H^{p,0}(X)\\
\bigwedge^pH^{1,0}(Z)= \!\!\!\!\!\!\!\!\!\!\!\!\!\!\!\!\!\!\!\!&H^{p,0}(\Alb(Z))\ar[r]^{\psi^Z_p}\ar[u]_{h^*}& H^{p,0}(Z)\ar[u]_{g^*}}
$$
\smallskip
Note that he vertical arrows are injective. From this diagram we see that if $p> \dim \alb(Z)$ then $\bigwedge^pH^{1,0}(Z)\subseteq \ker \psi_p$.
\end{remark}

\begin{remark}
When $h^{p,0}(Y)=0$ for all $p\geq 0$ - as happens for instance if $Y$ is rational - it is then possible to produce, via our method, elements in the kernel of $\psi_p $ for any $0\leq p\leq q$.
More generally, our method can be used to produce varieties $X$ with non-trivial elements in the kernel of the cup product map
$
H^{h,0}(X)\otimes H^{k,0}(X)\longrightarrow H^{h+k,0}(X).
$
\end{remark}

\medskip

We end this section by giving an explicit construction of $r={q \choose p}$ independent $p$-forms on $X$ contained in $\ker \psi_p$. The idea is to present $p$-forms invariant under the action of the monodromy group $S_d$. Let us fix a basis $\eta_1,\ldots,\eta_q\in H^{1,0}(Z)$ and let
\begin{displaymath}
\omega_{j}^{i}:=\pi_i^*(\eta_j)\in H^{1,0}(X)
\end{displaymath}
with $i=1,\ldots,d$ and $j=1,\ldots,q$. Then the following holds.
\begin{proposition}\label{proposition p-FORMS}
Under the assumption of Theorem \ref{theorem KERNEL}, for any choice of $p$ one-forms $\eta_{j_1},\ldots,\eta_{j_p}\in H^{1,0}(Z)$ of the basis,
\begin{equation}\label{equation LAGRANGIAN FORM}
\omega:=\sum_{i=1}^d \omega^{i}_{j_1}\wedge \omega^{i}_{j_2}\wedge \ldots \wedge \omega^{i}_{j_p}\in \bigwedge^p H^{1,0}(X)
\end{equation}
is an element of a copy of $A_p\subset \ker \psi_p$. Moreover, when $p=2$, the rank of  $\omega$  is $2(d-1)$.
\begin{proof}
As in the proof of Theorem \ref{theorem qCOPIE}, for any $j=1,\ldots, q$, we can define a copy $\Gamma_{\eta_j}\subset H^{1,0}(X)$ of the standard representation of $S_d$ associated to $\eta_j$, by pulling back $\eta_j$ via the $\pi_i$'s. So, let us consider the standard representations $\Gamma_{\eta_{j_1}},\ldots,\Gamma_{\eta_{j_p}}$ and the corresponding vector space $A_p$. Clearly, the $p$-form $\omega$ is invariant under the action of $S_d$ and it is contained in $A_p\in \ker \psi_p$.\\
Now, let $p=2$ and let $\omega=\sum_{i=1}^d \omega^{i}_{j_1}\wedge \omega^{i}_{j_2}\in \bigwedge^2 H^{1,0}(X)$. By the equations defining $\Gamma_{\eta_{j_1}}$ and $\Gamma_{\eta_{j_2}}$ we have
\begin{displaymath}
\omega^{d}_{j_1}=-\omega^{1}_{j_1}-\ldots -\omega^{d-1}_{j_1}\quad\textrm{and}\quad\omega^{d}_{j_2}=-\omega^{1}_{j_2}-\ldots -\omega^{d-1}_{j_2}.
\end{displaymath}
and hence
\begin{displaymath}
\omega=2\sum_{i=1}^{d-1} \omega^{i}_{j_1}\wedge \omega^{i}_{j_2}+\sum_{i\neq k} \omega^{i}_{j_1}\wedge \omega^{k}_{j_2}.
\end{displaymath}
To conclude, notice that $\omega^{1}_{j_1},\ldots,\omega^{d-1}_{j_1},\omega^{1}_{j_2},\ldots,\omega^{d-1}_{j_2}$ are independent forms on $X$ and the rank of $\omega$ is equal to the rank of the associated matrix of dimension $2(d-1)\times 2(d-1)$
\begin{displaymath}
A_{\omega}=\left[\begin{array}{cc}
0 & B_{\omega} \\
- B_{\omega} & 0
\end{array}\right],
\quad\quad\textrm{where}\quad
B_{\omega}=\left[\begin{array}{cccc}
1 & \frac{1}{2} & \ldots & \frac{1}{2} \\
\frac{1}{2} & 1 & \ldots & \frac{1}{2} \\
 &  & \ddots &  \\
\frac{1}{2} & \frac{1}{2} & \ldots & 1
\end{array}\right].
\end{displaymath}
\end{proof}
\end{proposition}
Note that when $p=2$, the form $\omega\in A_p$ is the element corresponding to the elementary symmetric function $-\xi_2$, under the identification (\ref{equation IDENTIFICATION}).

\bigskip
\section{Abelian surfaces of type $(1,2)$}\label{section BARTH}

Let us consider a couple   $(S,\mathcal L)$, where $S$ is a smooth complex abelian surface and $\mathcal L$ a line bundle over $S$ of degree $4$.
The line bundle $\mathcal L$ defines  a $(1,2)$-polarization on $S$, which we will denote $[\mathcal L]$.
In \cite{Barth}, Barth gives a detailed treatment of these surfaces, and  we here recall the results we need.

There is a natural isogeny associated to $\mathcal L$ - sometimes called itself the polarization - $\lambda_{\mathcal L} \colon S\longrightarrow \pic^0(S)=S^{\vee}$ defined by associating to $t\in S$ the invertible sheaf $\mathcal{L}^{-1}\otimes t^*\mathcal L$, where $t^*$ is the translation by $t$. We recall that \cite{Mumford}
\begin{equation}\label{tilambda}
T(\mathcal L):=\ker \lambda_{\mathcal L}\cong \mathbb Z/2\mathbb Z\times \mathbb Z/2\mathbb Z.
\end{equation}
It follows that $ T(\mathcal L)=\{x_0,x_1,x_2,x_3\}$, where $x_0$ is the origin of $S$, and $x_1, x_2, x_3$ are three points of order $2$.
As $h^0(S, \mathcal L)=\frac{1}{2} c_1(\mathcal L)^2=2$, the linear system $\vert \mathcal L\vert$ induces a linear pencil on $S$.

Let  $\mathcal{W}(1,2)$ the moduli space of abelian surfaces with a polarization of type $(1,2)$.
Given a couple $(S,[\mathcal L])\in \mathcal W(1,2)$, there exists an irreducible curve $C\in \vert \mathcal L\vert $ if and only if $\mathcal L$ is not of the form $\mathcal O_S(E+2F)$ where $E$ and $F$ are elliptic curves in $S$ such that $E^{\cdot 2}=F^{\cdot 2}=0$ and $E\cdot F=1$.
Moreover, if there exists an irreducible element in $\vert \mathcal L\vert $, then the general member is smooth, and the linear pencil $\vert \mathcal L\vert$ has precisely $T(\mathcal L)=\{x_0,x_1,x_2,x_3\}$ as  base points.

Let us suppose that there exists such an irreducible curve $C\in \vert \mathcal L\vert $.
In this case $C$ is either smooth of genus $3$, or an irreducible curve of geometric genus $2$ with one double point, which is easily seen to be a node.

Another important result in \cite{Barth} is that the  $(-1)$-involution on $S$ restricts to an involution on any curve $C\in \vert \mathcal L\vert$.
In the following we shall denote the induced involution  on $C$ as  $\iota$, and as $\pi$ the quotient morphism $C\longrightarrow C/\langle \iota\rangle$.

The surface  $S$ can be recovered from the data of the morphism $\pi$.  Indeed $S$ is the generalized Prym variety $\prym$ associated to this morphism (see \cite{L-B}, Section 12.3).  In order to fix the notation, we briefly define $\prym$ in our cases.
Let us distinguish the smooth and singular case.

\medskip

Suppose that  $C\in \vert \mathcal L\vert$ is a smooth curve.
By the Riemann-Hurwitz formula, the quotient $C/\langle \iota\rangle$ is a smooth elliptic curve $E$.
Consider the embedding of $E$ in the second symmetric product $C^{(2)}$
$$
E\cong \{p+\iota(p), p\in C\}\subset C^{(2)},
$$
and compose this map with the Abel map $C^{(2)}\hookrightarrow J(C)$. This is just  the inclusion given by pullback on the Picard varieties
$$
\pi^*\colon Pic^0(E)\cong E\hookrightarrow Pic^0(C)\cong J(C).
$$
Then, by composing the Jacobian embedding with the quotient map, we have a well defined morphism $\eta\colon C\longrightarrow J(C)/\pi^*E$, which is an embedding satisfying $\eta(C)^{\cdot 2}=4$  \cite[Proposition (1.8)]{Barth}.
The abelian surface $J(C)/\pi^* E$ is the generalized Prym variety $\prym$ associated to $\pi$; this notation is dual to the one used by Barth in section (1.4) of \cite{Barth}, see in particular the Duality Theorem (1.2).

\begin{proposition}[Barth]\label{proposition SMOOTH}
Let  $C\in \vert \mathcal L\vert$ a smooth curve.
Then there exists a degree 2 morphism to an elliptic curve $\pi\colon C\longrightarrow E$ such that $S$ is naturally identified with $\prym =J(C)/\pi^*E$.
Conversely, any smooth bielliptic  genus $3$ curve $\pi\colon C\longrightarrow E$ is embedded in  $\prym $ as a curve of self-intersection $4$.
\end{proposition}

\medskip

Let us consider the case when  $C\in \vert \mathcal L\vert$ is a singular irreducible curve.
As recalled above, $C$ has geometric genus $2$ with one node $q\in C$.
Let us call $q_1, q_2\in  \widetilde C$ the preimages of $q$.

The $(-1)$-involution on $S$ extends to an involution $\iota\colon C\longrightarrow C$ that fixes $x_0,\ldots x_3$.
Clearly, also $q$ is fixed by $\iota$ and the quotient $C/\langle \iota \rangle$ is an irreducible curve of arithmetic genus 1 with one node.
It is clear that the points $\{\nu^{-1}(x_0),\ldots, \nu^{-1}(x_3), q_1,q_2\}$ are the Weierstrass points of $\widetilde C$.

The isogeny $\varphi\colon J(\widetilde C)\longrightarrow S$ has degree $2$.
Indeed, the curve $\widetilde C$ has self-intersection $2$ in $J(\widetilde C)$, while $\nu_*\widetilde C^{\cdot 2}=C^{\cdot 2}=4$ in $S$.
Hence there exists a torsion point $\epsilon\in J(\widetilde C)$ such that $\varphi$ is the quotient map induced by the involution $z\mapsto z+\epsilon$.
Let us identify as usual $J(\widetilde C)$ with $\pic^0 (\widetilde C)$. Clearly $\epsilon \sim q_1-q_2$ in  $ \widetilde C$.
Indeed, $q_1-q_2$ is $2$-torsion because the $q_i$'s are Weierstrass points, and these two points are identified in $S$.

Let us call $G\cong \Z/2\Z $ the order $2$ subgroup of $J(\widetilde C)$ generated by $q_1-q_2$.
We shall call the quotient $J( \widetilde C)/G$ the generalized Prym variety $\prym$ associated to $\pi$.
We have then the following.

\begin{lemma}\label{singular}
Let $(S,\mathcal L)$  be such that there exists an irreducible singular curve $C\in \vert \mathcal L\vert$.
Then $S$ is naturally identified with  $\prym=J( \widetilde C)/G$.
\end{lemma}

\medskip
\section{Torsion points and geometry of the pencil $\vert \mathcal L\vert$}\label{section PENCIL}

Let  $(S,\mathcal L)$ be such that there exists an irreducible $C\in \vert \mathcal L\vert$.
As recalled above, the general element in $\vert \mathcal L\vert$ is smooth;
moreover, it is non-hyperelliptic (\cite[Claim 2.5]{TY}, \cite{PIETRO}).

We shall suppose  that any element of $\vert \mathcal L\vert$ is irreducible.
This amounts to asking that there are no curves of the form $E_1+E_2$, where $E_1$ and $E_2$ are smooth elliptic curves contained in $S$ meeting in two nodes.
This is a general condition in $\mathcal W(1,2)$.
Note that the above conditions would follow from the  requirement that  $S$  be a simple abelian surface, but this last assumption is much stronger, being generic in  $\mathcal W(1,2)$ \cite{PIETRO}.

\medskip

Under this assumption, we shall see that $\vert \mathcal L\vert$ has exactly 12 singular  members, corresponding to the order two points of $S$ different from the $x_i$'s. (Proposition \ref{proposition NODAL FIBERS}).
The linear system contains exactly $6$ smooth hyperelliptic elements that are related to a particular subset of the points of order $4$ (Lemma \ref{proposition SMOOTH} and Proposition \ref{proposition HYPERELLIPTIC FIBERS}).
In Proposition \ref{proposition ORDER THREE} we give a characterization of the triple points of $S$ in terms of the canonical images of the corresponding curves of the linear pencil.
\medskip

To fix the notation, let $\widetilde S$ be the blow up of $S$ in the base points $\{x_0,\ldots, x_3\}$ and
let $f\colon \widetilde S\longrightarrow \pr^1$ denote the fibration induced by the pencil.
We denote by $E_0, \ldots, E_3$ the four exceptional curves of the blow up, that are sections of $f$.

\begin{proposition}\label{proposition NODAL FIBERS}
Let  $(S,\mathcal L)$ be such that any element of $\vert \mathcal L\vert$ is irreducible.
The linear pencil  has $12$ singular elements that are all irreducible curves of  geometric genus $2$ with one node.
These nodes are the points of $S$ of order $2$ different from the $x_i$'s.
\begin{proof}
The first part of the proposition has already been established in Section \ref{section BARTH}.
We saw also that the singular points of the curves in $\vert \mathcal L\vert $ are points of order $2$ in $S$.
As $S$ has $16$ points of order two, $4$ of which are the base points $x_i$'s,
it remains only to prove is that there are $12$ singular curves in the linear pencil.
This derives from the following formula on the invariants of the fibration $f\colon \widetilde S\longrightarrow \pr^1$ \cite[Lemma VI.4]{Beauville}
$$
c_2({\widetilde S})=e(\pr^1)e(F)+\sum (e(N)-e(F)),
$$
where $F$ is a smooth fiber, and the sum is taken on all the singular fibers $N$ of $f$.
As any $N$ is irreducible with one node, we have that its topological characteristic is $e(N)=e(F)+1$.
Hence, if $n$ is the number of singular fibers, the formula becomes $4=-8+n$, and the proof is concluded.
\end{proof}
\end{proposition}

\smallskip

The other special elements of the linear pencil $\vert \mathcal L\vert$ are smooth hyperelliptic curves.
Let us consider the sets of points of $S$ defined as follows
$$
P_i:=\left\{ x\in S\, \,\, | \,\,\, 2x=x_i\right\}, \quad\quad \mbox{ for }i=1,2,3.
$$
Any of the $P_i$'s is a set of $16$ particular order $4$ points in $S$.
\begin{lemma}\label{smooth}
Let $x\in P_i$, and $D$ be the element of $\vert \mathcal L\vert $ passing through $x$.
Then $D$ is smooth.
\begin{proof}
Suppose by contradiction that $D$ is a singular curve.
Let $\nu\colon \widetilde D\longrightarrow D$ be its normalization, and $q_1,q_2\in \widetilde D$ be the inverse images of the node. By abuse of notation, let us still call $x_i$'s the inverse images $\nu^{-1}(x_i)$'s, and $x$ the inverse image of $x$ in $\widetilde D$.
Recall that $\{x_0,x_1,x_2, x_3, q_1,q_2\}$, are the Weierstrass points of the genus $2$ curve $\widetilde D$.
Recall from Lemma \ref{singular} that $S$ is naturally identified with $J(\widetilde D)/G$ where $G$ is the order two subgroup of $J(\widetilde D)$ generated by $q_1-q_2$.
The equality $2x=x_i$ in $S$, means in terms of this identification that
$$
2x\sim_{\widetilde D}x_0+x_i,\quad \mbox{ or }\quad
2x\sim_{\widetilde D}x_0+x_i+q_1-q_2.
$$
The first formula is impossible because $x_0+x_i$ (for $i\not =0$) does not belong to the $g^1_2$.
Applying the hyperelliptic involution  $\sigma$ to the second linear equivalence, we would obtain $2x\sim_{\widetilde D} 2\sigma (x)$.
This would imply that $x$ is a Weiestrass point of $\widetilde D$, which is a contradiction.
\end{proof}
\end{lemma}

\begin{proposition}\label{proposition HYPERELLIPTIC FIBERS}
Let  $(S,\mathcal L)$ be such that  any element of $\vert \mathcal L\vert$ is irreducible.
The linear pencil has  $6$ smooth hyperelliptic elements such that any $P_i$ consists of the Weierstrass points of two of them.
\begin{proof}
Let  $C$ be a smooth hyperelliptic curve belonging to $\vert \mathcal L\vert$.
By Barth's construction such a curve has a bielliptic involution $\iota\colon C\longrightarrow C$.
As  $\iota$ and $j$ commute,  $j$ induces a permutation on the fixed points of $\iota$, which are exactly the $x_i$'s that  does not fix any of them.
Let now $x\in C$ be a Weierstrass point. By what observed above, $2x\sim_C x_0+ j(x_0)=x_0+x_i$ for some $i\in\{1,2,3\}$.
Hence $x$ belongs to the set $P_i$.

On the other hand, let $x$ be a point in $P_1$, and $D\in \vert \mathcal L\vert$ be the curve passing through $x$.
We proved in Lemma \ref{proposition SMOOTH} that $D$ is smooth.
In particular we can identify $S$ with $J(D)/\pi^*E$, where $E$ is the quotient of $D$ by the bielliptic involution.
Using the above identification, we have that there exists $s\in D$ such that
$$2x-2x_0\sim_D x_1+s+\iota(s)-3x_0$$
 in $D$.
This implies that $2x\sim_D 2\iota(x)$; as $x$ is not fixed by the bielliptic involution $\iota$, $2x$ induces a $g^1_2$  on $D$.
Hence $D$ is hyperelliptic, and $x$ is one of its Weierstrass points.  Moreover, the other $7$ Weierstrass points necessarily satisfy $2y= x_1$ in $S$, so they also lie in $P_1$.

Choosing one point $x'\in P_1\setminus D\cap P_1$, we obtain another hyperelliptic curve $D'\in \vert \mathcal L\vert$, whose Weierstrass points are precisely $P_1\setminus D\cap P_1$.

Making the same construction for $x_2$ and  $x_3$, we obtain the other $4$ hyperelliptic curves, and the statement is verified.
\end{proof}
\end{proposition}

\begin{remark}
Recall that $f\colon \widetilde S\longrightarrow \pr^1$  is the fibration obtained blowing up the base points of the pencil $\vert \mathcal L\vert$. Consider the morphism of sheaves
$$
\sym^2 f_*\om_{f}\stackrel{\delta}{\longrightarrow} f_*\om_{f}^2.
$$
The fiber of $\delta $ over $p\in \pr^1$ is  the multiplication morphism
$$
\sym^2 H^0(f^{-1}(p), \om_{f^{-1}(p)})\longrightarrow H^0(f^{-1}(p), \om _{f^{-1}(p)}^2).
$$

The cokernel of $\delta$ is a torsion sheaf  $\mathcal T$ supported over the points corresponding to the smooth hyperelliptic fibers of $f$.
As Reid proves in \cite[Sections 3.2 and 3.3]{Reid}, the following relation holds
$$
K^2_{\widetilde S}=3\chi(\oo_{\widetilde S}) -10  +\deg \mathcal T,
$$
hence in our case the degree of $ \mathcal T$ is $6$. In \cite{Reid} it is also proved that  the contribution of any hyperelliptic smooth fiber to the degree of $\mathcal T$, which is usually called Horikawa number, is greater or equal to one. We can therefore conclude that these fibers have Horikawa numbers equal to $1$, i.e. that  the modular image of the base $\pr^1$ in  $\overline{\mathcal M}_3$ intersects the hyperelliptic locus transversally in the points corresponding to these fibers.
\end{remark}

\medskip

With the following result we give a geometric description - in terms of the linear pencil - of the points of order 3 of $S$.
Before stating the result, let us fix a notation. Let $C$ be an irreducible curve of arithmetic genus $3$ with one node $q$. Let $\nu\colon \widetilde C\longrightarrow C$ be the normalization, and $q_1, q_2 $ be the preimages of the node.
The  canonical embedding of $C$ in $\pr^2$ is associated to its dualizing sheaf $\nu_*(\omega_{\widetilde C}(q_1+q_2))$, and corresponds on $\widetilde C$ to the birational morphism induced by the divisor $K_{\widetilde C}+q_1+q_2$.

\begin{proposition}\label{proposition ORDER THREE}
Let $p\neq x_0$ be a point of $S$ and let $C\in \vert \mathcal L \vert$ be the curve passing through $p$. The following are equivalent:
\begin{itemize}
  \item[(i)] $p$ is a point of order 3 of $S$;
  \item[(ii)] $C$ is not a smooth and hyperelliptic  curve and, when  its  image via the canonical embedding  $\varphi\colon C\hookrightarrow  \pr^2$ has an inflection point of order 3 at $\varphi(p)$ with tangent line $\overline{\varphi(p_0)\varphi(p)}$.
\end{itemize}
\begin{proof}
$(i)\Rightarrow (ii)$.
Suppose that $C$ is the element of the pencil passing through a point $p$ of order three.
If $C$ is smooth, let $E=C/\langle \iota \rangle$ be its bielliptic quotient. Using the identification $S=J(C)/\pi^*E$,  the assumption implies that there exists a point $s\in C$ such that
$$3(p-x_0)\sim_C s+\iota(s)-2x_0,$$
that is $3p\sim_C s+\iota (s) +x_0$ and this  induces a base point free $g_3^1$ on $C$. So, in particular, assumption (i) implies that $C$ is not a smooth hyperelliptic curve.
Let us distinguish between the smooth and the singular case.

Assume that $C$ is  smooth and non-hyperelliptic.
As observed above, there exists a point $s\in C$ such that  $3p\sim_C x_0+s+\iota(s)$.
As $C$ is non-hyperelliptic, there exists $r\in C$ such that
$$
3p+r \sim_C x_0+s+\iota(s)+r\sim_C K_C.
$$
The latter equivalence proves that $\varphi(p)$ is an inflection points of order 3 for the canonical image of $C$ in $\pr^2$ with tangent line $\overline{\varphi(r)\varphi(p)}$. To complete the first part of the proof we need to show that $r=x_0$. Observe that
$$
x_0+\iota(s)+s+r\sim_C K_C\sim \iota K_C\sim_C x_0+s+\iota(s)+\iota(r).
$$
So, $r=\iota(r)$ and hence $r$ is one of the $x_i$'s. On the other hand $K_C\sim_C x_0+x_1+x_2+x_3$, then the only possibility is $r=x_0$.

Let us now suppose that $C$ is nodal. As usual, we identify $S$ with $J(\widetilde C)/G$.
As $p$ has order $3$ in $S$, we have that $3p\sim_{\widetilde C}3x_0+q_1-q_2$ in $\widetilde C$.
Now note that $q_2\sim_{\widetilde C}2x_0-q_2$, because both $x_0$ and $q_2$ are Weierstrass points.
Hence,
\begin{equation}\label{eh}
3p+x_0\sim_{\widetilde C}2x_0+q_1+q_2\sim_{\widetilde C} K_{\widetilde C}+q_1+q_2,
\end{equation}
as wanted.

$(ii)\Rightarrow (i)$. Suppose $C$ is not simultaneously smooth and hyperelliptic, and that  its canonical image has an inflection point of order 3 at $\varphi(p)$ with tangent line $\overline{\varphi(x_0)\varphi(p)}$.

If $C$ is smooth, we have that $3p+x_0\sim_C K_C$.
For some points $a,b\in C$, let $2x_0+a+b\sim_C K_C$ be the divisor cut out by the tangent line to $C$ at $x_0$.
Notice that $2x_0+a+b\sim_C K_C\sim_C \iota K_C\sim_C 2x_0+\iota(a)+\iota(b)$.
Therefore either $a=\iota(b)$ with $a\neq b$ or $a=\iota(a)$.
In any case, the relation $3(p-x_0)\sim_C s+\iota(s)-2x_0$ holds for some $s\in C$, so $3p=0$ in $S$.

In the singular case, we have $3p+x_0\sim_{\widetilde C}K_{\widetilde C}+q_1+q_2$, which implies the statement by equation (\ref{eh}).
\end{proof}
\end{proposition}

\medskip


\section{The triple covering construction}\label{section TRIPLE}

Recall that $\widetilde S$ is the blow up of $S$ in $\{x_0,x_1,x_2,x_3\}$, and $f\colon \widetilde S\longrightarrow \pr^1$ the induced fibration. We denote by $E_0, \ldots, E_3$ the exceptional divisors, which are sections of $f$.
Given a fiber $F$ of $f$ we shall always denote by $p_i$ the point of intersection between $E_i$ and $F$; so, when we identify the fiber $F$ with the corresponding element of the linear pencil $\vert \mathcal L\vert$, the $p_i$'s correspond to the base points $x_i$'s.

The homomorphism of sheaves ${f}^*{f}_*\om_{f}(-E_0)\longrightarrow \om_{f}(-E_0),$
induces a  relative rational map
$$
\xymatrix{ \widetilde S \ar@{-->}[rr]^-{\gamma} \ar[dr]^{f} & &Y \ar[dl] \\ & \pr^1 & \\ }
$$
Where $Y:=\pr (f_*\om_{f}(-E_0))$ is the relative projective bundle on $\pr^1$ associated to the rank $2$ vector bundle $f_*\om_{f}(-E_0)$. Hence $Y$ is a rational surface.

We note that $\gamma$ is a generically finite 3 to 1 map. Indeed, on a smooth non-hyperelliptic fiber $F$ of $f$, the restriction $\gamma_{|F}\colon F \longrightarrow \mathbb{P}^1$ corresponds to the composition of the canonical immersion $\varphi\colon C\hookrightarrow \pr^2$  with the projection from $\varphi(p_0)$.

\begin{proposition}\label{lemma HIRZEBRUCH}
With the same notation as above, the vector bundle ${f}_*\om_{f}(-E_0)$ is isomorphic to  $\oo_{\pr^1}(2)\oplus \oo_{\pr^1}(-1)$.
Hence  $Y$ is isomorphic to the minimal rational surface $\mathbb F_3=\pr(\oo_{\pr^1}(3)\oplus \oo_{\pr^1})$
\begin{proof}
Observe firstly that we have the following decomposition (see for instance \cite{Fujita}).
$$
{f}_*\om_{f}= \oo_{\pr^1}(\alpha) \oplus \oo_{\pr^1}^{\oplus 2},
$$
where $\alpha = \deg {f}_*\om_{f}= \chi(\oo_{\widetilde S})-\chi(\oo_F)\chi(\oo_{\pr^1})=2$.
By Grauert's Theorem \cite[Chapter III Corollary 12.9]{Hartshorne} the sheaf $R^1{f}_*\oo_{\widetilde S}(E_0)$ is locally free, hence by relative duality
$$
R^1{f}_*\oo_{\widetilde S}(E_0)\cong \left({f}_*\om_{f}(-E_0)\right)^{\vee}.
$$
Let us consider the short exact sequence of sheaves
$$
0\longrightarrow \oo_{\widetilde S}\longrightarrow \oo_{\widetilde S}(E_0)\longrightarrow \oo_{E_0}(E_0)\longrightarrow 0
$$
and the long exact sequence induced by the pushforward
\begin{equation}\label{long sequence}
0\longrightarrow {f}_* \oo_{\widetilde S}\longrightarrow {f}_* \oo_{\widetilde S}(E_0)\longrightarrow {f}_*\oo_{E_0}(E_0) \longrightarrow  R^1{f}_*\oo_{\widetilde S}\longrightarrow  R^1{f}_*\oo_{\widetilde S}(E_0)\longrightarrow 0.
\end{equation}
Observe that ${f}_* \oo_{\widetilde S}\cong \oo_{\pr^1}$ and that ${f}_*\oo_{E_0}(E_0) \cong \oo_{\pr^1}(-1)$.
Moreover, by relative duality again, we have that
$$
R^1{f}_*\oo_{\widetilde S}\cong \left( {f}_*\om_{f}\right)^{\vee}= \oo_{\pr^1}\oplus \oo_{\pr^1}\oplus \oo_{\pr^1}(-2).
$$
Hence, from (\ref{long sequence}) we deduce the exact sequence
$$
0\longrightarrow \oo_{\pr^1}(-1) \longrightarrow\oo_{\pr^1}\oplus \oo_{\pr^1}\oplus \oo_{\pr^1}(-2) \longrightarrow R^1{f}_*\oo_{\widetilde S}(E_0)\longrightarrow 0.
$$

We deduce that the image of the morphism $\oo_{\pr^1}(-1) \longrightarrow\oo_{\pr^1}\oplus \oo_{\pr^1}\oplus \oo_{\pr^1}(-2)$ is contained in  $\oo_{\pr^1}\oplus \oo_{\pr^1}$.
Hence $ \oo_{\pr^1}(-2)$ injects into $R^1{f}_*\oo_{\widetilde S}(E_0)$ and $ R^1{f}_*\oo_{\widetilde S}(E_0)= \oo_{\pr^1}(-2)\oplus  \oo_{\pr^1}(\beta)$ for some $\beta$. Finally, by computing the degrees of these sheaves, we see that $\beta= 1$.
\end{proof}
\end{proposition}

The rational map $\gamma\colon \widetilde S\dashrightarrow Y $ is not a morphism.
Indeed, a point $b\in \widetilde S$ is an indeterminacy point for $\gamma$ if and only if the associated morphism of sheaves ${f}^*{f}_*\om_{f}(-E_0)\longrightarrow \om_{f}(-E_0)$ is not surjective at $b$  \cite[Chapter II, Section 7]{Hartshorne}.
By Nakayama's lemma, the morphism is surjective if and only if its restrictions to the fiber of $f$ are. For a fiber $F$ this restriction is the evaluation morphism
$$
H^0(F,\om_F(-p_0))\otimes \oo_F \longrightarrow \om_F(-p_0),
$$
which is surjective if and only if the line bundle $\om_F(-p_0)$ is globally generated. On the smooth non-hyperelliptic fibers, as well as on the singular ones, it is easy to check that this is the case.
Let $D$ be a smooth hyperelliptic fiber.
By Proposition \ref{proposition HYPERELLIPTIC FIBERS} we have that the hyperelliptic involution maps $p_0$ in one of the $p_k$ for some $1\leq k\leq 3$, and $p_k$ is a base point for the linear system $\vert \om_{D}(-p_0)\vert \cong\om_{f}(-E_0)_{|{D}}$.
Recalling that by Proposition \ref{proposition HYPERELLIPTIC FIBERS} there are $6$ hyperelliptic fibers of $f$,  we have proven the following

\begin{proposition}
The rational map $\gamma$ has $6$ indeterminacy points, one on any hyperelliptic fiber of $f$.
Moreover, any section $E_k$ contains two of them.
\end{proposition}

\begin{notation}\label{notation ON barS}
For $j=1,2$ and $1\leq k\leq 3$, let us denote each of these indeterminacy points by $b_{jk}$, with the convention that $b_{1k}, b_{2k}\in E_k$.
Let
$$
\xymatrix{ \overline{S} \ar[rr] \ar[dr]_{\bar{f}} & & \widetilde S \ar[dl]^f \\ & \mathbb{P}^1 & \\ }
$$
be the blow up of $\widetilde{S}$ at the $b_{jk}$'s, together with the induced fibration $\bar{f}$.
We denote by $G_{jk}$ be the exceptional divisor corresponding to the point $b_{jk}$. Moreover, by abuse of notation let $E_0,\ldots,E_3\subset \overline{S}$ be the strict transforms of the sections $E_0,\ldots,E_3\subset \widetilde{S}$ and let $D\subset\overline{S}$ be the strict transform of any hyperelliptic fiber $D\subset\widetilde{S}$.
\end{notation}

\begin{remark}\label{fibre}
Summing up, we have proved above that any fiber of $\overline{f}\colon \overline{S}\longrightarrow Y$ is one of the following:
\begin{itemize}
  \item[-] a smooth irreducible curve of genus 3;
  \item[-] a nodal irreducible curve of geometric genus 2 as in Proposition \ref{proposition NODAL FIBERS};
  \item[-] a reducible nodal curve $F=D\cup G_{jk}$, where $D$ is a smooth hyperelliptic curve and $G_{jk}\cong \mathbb{P}^1$ (see Proposition \ref{proposition HYPERELLIPTIC FIBERS} and Notation \ref{notation ON barS}).
\end{itemize}
\end{remark}

\begin{proposition}\label{proposition BLOW UP}
The sheaf $\bar{f}_*\oo_{\overline{S}}(E_1+E_2+E_3)$ induces a finite degree $3$ morphism $\overline \gamma$
$$
\xymatrix{ \overline S \ar[rr]^-{\overline{\gamma}} \ar[dr]_{\bar{f}} & & Y \ar[dl] \\ & \mathbb{P}^1 & \\ }
$$
resolving the indeterminacy of $\gamma$.
\begin{proof}
Recall that $\gamma$ is induced by the sheaf $f_*\om_f(-E_0)=f_*\oo_{\widetilde{S}}(E_1+E_2+E_3)$. We just have to show that $\overline{\gamma}$ restricts to a well defined morphism on any fiber. Away from the $G_{jk}$'s $\overline{\gamma}$ coincides with $\gamma$.
Let $G=G_{jk}$ and let $D$ be the hyperelliptic fiber passing through $b_{jk}$.
Without loss of generality, let $k=1$.
Let us consider the total transform of $D\subset\widetilde{S}$, which is given by $D\cup G\subset \overline S$.
The sheaf defining the restriction of $\overline{\gamma}$ to $G$ is $\oo_{\overline{S}}(E_1+E_2+E_3)_{|G}\cong \oo_{G}(1)$, and hence
 $\overline{\gamma}_{|{G}}\colon G\longrightarrow\mathbb{P}^1$ is an isomorphism.
On the other hand, the restriction of the map $\overline{\gamma}$ to $D$ is given by the sheaf
$$
\oo_{\overline{S}}(E_1+E_2+E_3)_{|D}=\oo_D\left( (E_2\cap D)+(E_3\cap D) \right)=\oo_D(p_2+p_3).
$$
By Proposition \ref{proposition HYPERELLIPTIC FIBERS}, $p_2$ and $p_3$ are conjugate under the hyperelliptic involution of $D$. Hence the linear system $|p_2+p_3|$ on $D$ is the $g^1_2$. Therefore $\overline{\gamma}_{|D\cup G}$ has no base points and turns out to be a degree three morphism to $\mathbb{P}^1$, as in Figure \ref{figure HYPERELLIPTIC FIBER} below, where $p_0$ denotes the intersection of the fiber with $E_0$, $b_{jk}= G_{jk}\cap E_k$ and $\bar{b}_{jk}= G_{jk}\cap D$. To conclude, as $\overline \gamma$ does not contract any curve, it is a finite morphism of degree $3$.
\begin{center}
\begin{figure}[h!]
\input FIBRAiperellittica.pstex_t
\caption{}
\label{figure HYPERELLIPTIC FIBER}
\end{figure}
\end{center}
\end{proof}
\end{proposition}

We now compute the numerical equivalence class of the ramification divisor of $\overline \gamma$.

\begin{proposition}\label{lemma Rgamma}
The ramification divisor $\rg\subset \overline S$ is numerically equivalent to the divisor
$$
E_0+3\sum_{k=1}^3E_k+2\sum_{k=1}^3(G_{1k}+G_{2k})+5F,
$$
where $F$ is a fiber of $\overline f$.
\begin{proof}
The N\'eron-Severi group of $Y$ is generated by the  class of a fiber $\Gamma$ and by the class of the section with minimal self-intersection $C_0$. Moreover,
$K_Y\equiv -5\Gamma-2C_0$. By the formula for blow ups, we have
\begin{equation}\label{equation CANONICAL DIVISOR}
K_{\overline S}\equiv \sum_{i=0}^3E_i+2\sum_{k=1}^3(G_{1k}+G_{2k}).
\end{equation}
We remark that $\overline \gamma^*C_0=E_1+E_2+E_3$. To see this, notice that the sections $E_k$'s for $k\ne 0$ are $-3$-curves in $\overline S$. Moreover, they do not intersect the ramification locus $\rg$, because the $E_k$'s are disjoint curves and for any fiber $F$, the points $p_k=E_k\cap F$ map on the same point of $Y$. As the images of the $E_k$'s are $-3$-curves in $Y$, they are $C_0$.  Finally, $K_{\overline S}\equiv \rg +{\overline\gamma}^*K_Y$ by Riemann-Hurwitz formula, and the assertion follows.
\end{proof}
\end{proposition}

\medskip
\section{The Galois closure of $\overline \gamma$}\label{section GALOIS}

In this section we deal with the Galois closure of the degree three covering $\overline{\gamma}\colon \overline{S}\longrightarrow Y$ (see Definition \ref{definition GALOIS CLOSURE}). We note first that the monodromy group $M(\gamma)$ of the covering is the full symmetric group $S_3$. Indeed, if $M(\gamma)$ were isomorphic to $\mathbb{Z}/3\mathbb{Z}$, the covering would have only total ramification points, but this is not the case (see for instance Figure \ref{figure HYPERELLIPTIC FIBER}). Therefore the Galois closure of $\overline{\gamma}$ is the normalization of the Zariski closure $W$ of the variety (cf. Remark \ref{remark Z^(d-1)})
\begin{equation}\label{equation W0}
W^\circ:= \left\{ (p,q)\in \overline S\times \overline S \mid p\neq q\;\textrm{and}\;\overline \gamma (p)=\overline \gamma (q)\right\}\subset \overline S\times \overline S.
\end{equation}
In other words, the Zariski closure of $W^\circ$ is the divisor of the fibred product $\overline S\times_Y \overline S$ residual to the diagonal $\Delta$ of $\overline{S}\times \overline{S}$.
In the following we shall prove that the Zariski closure $W$ of $W^\circ$ is a normal surface itself and hence it is the Galois closure of $\overline{\gamma}$ (cf Proposition \ref{proposition Balpha SINGULARITIES}).

Let us denote by $\alpha_1,\alpha_2\colon W\subset \overline{S}\times \overline{S}\longrightarrow \overline{S}$ the natural projection maps and let $\alpha_3\colon W\longrightarrow \overline{S}$ be the morphism sending a point $(p,q)\in W$ to the point $r\in \overline{S}$ such that $\overline{\gamma}^{-1}(\overline\gamma (p))=\{p,q,r\}$. We note that the $\alpha_i$'s are generically finite morphisms of degree $2$. Given a general point $p\in \overline S$, the inverse image of $\overline \gamma(p)$ consists of three distinct points $\{p,q,r\}$ and hence $\alpha_1^{-1}(p)=\{(p,q),(p,r)\}$. Furthermore, since the $\alpha_i$'s do not contract any curve on $W$, we conclude that they are double coverings.
Now, let us consider the following commutative diagram
$$
\xymatrix{ W\ar[dr]\ar[d]_{\alpha_1} & \\ \overline S\ar[d]_{\overline f}\ar[r]^{\overline \gamma} & Y\ar[dl]\\ \pr^1 & \\}
$$
and let $\ba\subset \overline{S}$ denote the branch divisor of $\alpha_1$. We shall compute its numerical equivalence class and we shall prove that $\ba$ is a reduced curve with at most simple singularities.

Let $F\subset \overline{S}$ be a fiber of the morphism $\bar{f}:\overline{S}\longrightarrow \mathbb{P}^1$ and consider the restriction $\overline{\gamma}_{|F}\colon F\longrightarrow \mathbb{P}^1$ of the morphism $\overline{\gamma}$ to $F$. As usual, for $i=0,\ldots,3$, let $p_i=E_i\cap F$. We recall that when $F$ is irreducible, the map $\overline{\gamma}_{|F}$ is the projection of the canonical image of $F$ in $\pr^2$ from the point corresponding to $p_0$. For the sake of simplicity, hereafter we identify $F$ and its canonical image in $\mathbb{P}^2$. Let us define the subsets of $F$
\begin{equation*}
A:=\left\{p\in F\smallsetminus\{p_0\}\,|\, \exists\, q\in F\smallsetminus\{p,p_0\}\,:\,\overline{p_0p}\,\textrm{ is tangent at }q\right\}
\end{equation*}
and
\begin{equation*}
B:=\left\{p\in F\smallsetminus\{p_0\}\,|\,p\textrm{ is an inflection point of order 3 with tangent line }\overline{p_0p}\right\}.
\end{equation*}
The points of $A$ and $B$ correspond to the configurations (a) and (b) in Figure \ref{figure BRANCH POINTS} below.

On the other hand, let $F=D\cup G_{jk}$ where $D$ is a smooth hyperelliptic curve of genus 3. Then the restriction of $\overline{\gamma}$ to $F$ is described in Proposition \ref{proposition BLOW UP} (see also Figure \ref{figure HYPERELLIPTIC FIBER}). In particular, $\overline{\gamma}_{|D}\colon D\longrightarrow \mathbb{P}^1$ is the hyperelliptic map and $\overline{\gamma}_{|G_{jk}}\colon G_{jk}\longrightarrow \mathbb{P}^1$ is an isomorphism. So, we define the subset $T$ of $F$ given by
\begin{equation*}
T:=\left\{p\in G_{j,k}\,|\,\exists\, q\in D\,:\,q\textrm{ is a Weierstrass point and }\overline{\gamma}_{|F}(p)=\overline{\gamma}_{|F}(q)\right\}.
\end{equation*}

The following proposition describes the intersection of the branch curve with each fiber.

\begin{proposition}\label{proposition CONFIGURATIONS}
Let $F\subset \overline{S}$ be a fiber of the morphism $\bar{f}:\overline{S}\longrightarrow \mathbb{P}^1$. With the above notation, the intersection divisor $B_{\alpha_1|F}$ induced on $F$ by the branch divisor $\ba$ is given by one of the following.
\begin{itemize}
  \item[(i)] If $F$ is a smooth fiber, then $B_{\alpha_1|F}=\sum_{p\in A} p+2\sum_{p\in B}p+2p_0$ when $p_0=F\cap E_0$ is an inflection point of order 4 on $F$, and $B_{\alpha_1|F}=\sum_{p\in A} p+2\sum_{p\in B}p$ otherwise.
  \item[(ii)] If $F=N$ is a nodal irreducible fiber, then $B_{\alpha_1|F}=\sum_{p\in A} p+2\sum_{p\in B}p+2p_0$. In particular, the tangent line to $F$ at $p_0$ meets the node transversally.
  \item[(iii)] If $F=D\cup G_{jk}$ with $D$ hyperelliptic, then $B_{\alpha_1|F}=\sum_{p\in T} p+ 2p_0$.
\end{itemize}
\end{proposition}
Before proving the above result, let us state the following preliminary lemma.
\begin{lemma}\label{lemma MULTIPLICITY}
Let $p\in F\cap \ba$ for some fiber $F$ and let $(p,q)=\alpha_1^{-1}(p)\in W$. Then the multiplicity $m_p(\ba_{|F})$ of $\ba_{|F}$ at $p$ is equal to the multiplicity $m_q(\rg_{|F})$ of $\rg_{|F}$ at $q$.
\begin{proof}
To prove the assertion, it suffices to give a local description of $\ba$ in a neighborhood of a total ramification point $p\in \rg\cap \ba$ of $\overline{\gamma}$. So, let $(x,t)$ be local coordinates of $\overline{S}$ centered at $p$, where $t=constant$ is the local equation of a fiber of $\overline{f}$, and $x$ defines the coordinate along the fibers. The morphism $\overline{\gamma}$ is locally given by $(x,t)\mapsto \left(x^3+tr_1(t)x^2+tr_2(t)x+r_3(t),t\right)$, where $r_1,r_2,r_3\in\mathbb{C}[t]$. By (\ref{equation W0}) it is immediate to check that $W^\circ\subset \overline{S}\times \overline{S}$ is locally described by
$$\left\{ \left((x,t),(y,s)\right)\in \mathbb{C}^2\times \mathbb{C}^2\vert t=s, x^2+xy+y^2+t\left(p_1(t)(x+y)+p_2(t)\right)=0\right\}.$$
Therefore a point $(x,t)$ lies on $\ba$ if and only if the latter equation admits a unique solution with respect to $y$. Thus the local equation of $\ba\subset \overline{S}$ is $3x^2+2tp_1(t)x-t^2(p_1(t))^2+4tp_2(t)=0$. We conclude that around a total ramification point $p\in \ba$, the multiplicity of $\ba$ at $p$ is $2$, whereas it is $1$ at any other point of the neighborhood.
\end{proof}
\end{lemma}

\begin{proof}[\textit{Proof of Proposition \ref{proposition CONFIGURATIONS}}]
We study separately the three situations listed in the proposition.

\stepempty{i} Let $F$ be a smooth fiber and let $\overline{\gamma}_{|F}\colon F\longrightarrow \mathbb{P}^1$ be the projection of the canonical image of $F\subset \pr^2$ from the point corresponding to $p_0$. Let $p\in\ba\cap F$ for some fiber $F$ and let $(p,q)=\alpha_1^{-1}(p)\in W$. Then the canonical divisor $K_p\in Div\,F$ cut out by the line $\overline{p_0p}$ is one of the following (cf. Figure \ref{figure BRANCH POINTS} below):
\begin{center}
(a) $p_0+p+2q$\quad (b) $p_0+3p$\quad (c) $4p_0$\quad (d) $2p_0+2q$\quad (e) $3p_0+p$
\end{center}
\begin{center}
\begin{figure}[h!]
\input FIBREbranch.pstex_t
\caption{}
\label{figure BRANCH POINTS}
\end{figure}
\end{center}

Firstly we show that cases (d) and (e) cannot occur. Recall that the points $p_0,\ldots,p_3$ are collinear and they are the only fixed points under the action of the bielliptic involution $\iota$ on $F$. Moreover, such an involution is induced by a projectivity. So, if $\overline{p_0p}$ were a bi-tangent line as in (d) - that is $K_p=2p_0+2q$ with $p=p_0$ and $q\neq p_0$ - we would have that $\iota_*K_p=2p_0+2\iota(q)$. As the tangent line to $F$ at $p_0$ is tangent at $q$, we deduce that $q$ is fixed by $\iota$. Hence $q=p_i$ for some $i\neq 0$, but this is impossible because the $p_i$'s are collinear.

Analogously, suppose that $p_0$ is a flex point for $F$ as in case (e), that is $K_p=3p_0+p$ with $q=p_0$ and $p\neq p_0$. Thus $\iota_*K_p=3p_0+\iota(p)$ and $p=\iota(p)=p_i$ for some $i\neq 0$, a contradiction.

When $K_p=p_0+p+2q$ as in (a), then $p\in A$ and $q$ is a ramification point of index 2. Hence $m_p(\ba_{|F})=1$.

Case (b) happens if and only if $F$ has an inflection point of order 3 at $p$ with tangent line $\overline{p_0p}$, that is $p\in B$. In this case $p$ is a total ramification point of $\overline{\gamma}$ and hence $m_p(\ba_{|F})=2$ by Lemma \ref{lemma MULTIPLICITY}.


\medskip
\stepempty{ii} Let $F=N\subset \overline{S}$ be a nodal irreducible fiber and let $p\in N\cap \ba$. As the restriction $\overline{\gamma}_{|F}\colon F\longrightarrow \mathbb{P}^1$ is the projection from $p_0\in F\cap E_0$, we deduce that away from the node, two configurations analogous to (a) and (b) above are still possible. On the other hand, the cases (c), (d) and (e) cannot occur, whereas there is the following additional case: the tangent line at $p_0$ meets the node transversally (see (f) in Figure \ref{figure BRANCH POINTS}).

To see this fact, denote by $q$ the node of $N$ and let $L=\overline{p_0q}$ be the line through $p_0$ and $q$. We recall from Proposition \ref{proposition NODAL FIBERS} that $q\neq p_i$ for all $i$. Consider the points $\{p_0,p,q\}=L\cap N$. As $p_0$ and $q$ are fixed by the action of the involution $\iota$ on $N$, we have that also $p$ must be. Hence either $p=q$ or $p=p_0$. The same argument works for the line $L_i$ through $q$ and $p_i$, with $i=1,2,3$. We now prove that $p=p_0$.
Suppose by contradiction that $p=q$. Hence $L$ is one of the two tangent lines at $q$. Notice that the four lines $L,L_1,L_2,L_3$ must be distinct because $p_0,\ldots,p_3$ are collinear. Thus there exist two of those lines that are tangent in the corresponding $p_i$. Without loss of generality, let $L_1$ and $L_2$ be these lines. Let $\nu:\widetilde{\mathbb{P}}^2\longrightarrow \mathbb{P}^2$ be the blow up of $\mathbb{P}^2$ at $q$ and let $\widetilde{N}\subset \widetilde{\mathbb{P}}^2$ be the strict transform of $N\subset\mathbb{P}^2$. As in Proposition \ref{proposition NODAL FIBERS}, let $r_i=\nu^{-1}(p_i)$ and $\{q_1,q_2\}=\nu^*q$. Therefore $\nu^*(L\cap N)=r_0+q_1+2q_2$, $\nu^*(L_1\cap N)=2r_1+q_1+q_2$ and $\nu^*(L_2\cap N)=2r_2+q_1+q_2$ are linear equivalent divisors on $\widetilde{N}$. Thus $r_0+q_1$ is equivalent to $2r_1$, but this is impossible because $|2r_1|$ is the $g_2^1$ on $\widetilde{N}$ (see Proposition \ref{proposition NODAL FIBERS}).

Then we conclude that $p=p_0$, which is equivalent to configuration (f). Hence the section $E_0$ meet $N$ transversally at $p_0$ and the node $q\in \rg$ with $m_p(\ba_{|N})=2$. Moreover, this implies that cases (c), (d) and (e) are not possible. 

\medskip
\stepempty{iii} Let $F=D\cup G_{jk}$, where $D$ is a smooth hyperelliptic curve of genus 3. Let $w_1,\ldots,w_8\in D$ be the Weierstrass point and let $g_1,\ldots,g_8\in G_{jk}$ such that $\overline{\gamma}(w_t)=\overline{\gamma}(g_t)$ for any $t$. Hence the $g_t$'s lie on $\ba$ and for any $t$ we have that $w_t\in \rg$ with $m_{w_t}(\rg_{|F})=1$.

The last branch point on this fiber is the point $p_0$. Indeed, the hyperelliptic involution maps $p_0$ into the point $\bar{b}_{jk}:=D\cap G_{jk}$ (cf. Proposition \ref{proposition HYPERELLIPTIC FIBERS} and Figure 1). Hence $\bar{b}_{jk}\in\rg$ and it is a singular point of $F$. Thus $m_{\bar{b}_{jk}}(\rg_{|F})=2$ and  we are done. 
\end{proof}

We can now compute the numerical equivalence class of the branch locus $\ba$. We use  Notation \ref{notation ON barS}.
\begin{proposition}\label{proposition Balpha}
The branch divisor $\ba\subset \overline{S}$ is numerically equivalent to the divisor
\begin{equation}\label{equation CLASS OF Balpha}
-2E_0+4\sum_{k=1}^3E_k+20F-4\sum_{k=1}^3(G_{1k}+G_{2k}).
\end{equation}
\begin{proof}
Let $F$ be a general fiber of $\bar{f}\colon\overline{S}\longrightarrow \mathbb{P}^1$. We prove first that the restriction of $\ba+2\rg$ to a fiber is numerically equivalent to  $10(E_1+E_2+E_3)_{|F}$. The curve $F$ does not possess flex points by Proposition \ref{proposition ORDER THREE}. Hence the ramification divisor $\rg\subset\overline{S}$ meets $F$ at 10 points of ramification index 2. Let $q\in \rg\cap F$ and $p\in\ba\cap F$ such that $p_0+p+2q$ is the canonical divisor on $F$ cut out by the line $L=\overline{p_0q}$. As $(\ba+2\rg)_{|F\cap L}=p+2q\in |K_{F}(-p_0)|$, we have that $(\ba+2\rg)_{|F}$ is 10 times the $g^1_3$ defining $\overline{\gamma}_{|F}$, that is $(\ba+2\rg)_{|F}\equiv 10(E_1+E_2+E_3)_{|F}$.

Since $F\cdot G_{jk}=F\cdot F=0$ for any $j=1,2$ and $k=1,2,3$, there exist some integers $m,n_{jk}$ such that
\begin{align*}
\ba & \equiv -2\rg + 10(E_1+E_2+E_3) + mF + \sum_{k=1}^3(n_{1k}G_{1k}+n_{2k}G_{2k})\\
 & \equiv -2E_0 +4\sum_{k=1}^3E_k + (m-10)F + \sum_{k=1}^3((n_{1k}-4)G_{1k}+(n_{2k}-4)G_{2k}).
\end{align*}
Then $\ba\cdot F=10$ and from the description of Proposition \ref{proposition CONFIGURATIONS}, we have that $\ba\cdot G_{jk}=8$ and $\ba\cdot E_{k}=0$. Thus we deduce $m=30$ and $n_{jk}=0$ for any $j$ and $k$.
\end{proof}
\end{proposition}

\begin{proposition}\label{proposition Balpha SINGULARITIES}
The branch divisor $\ba$ is reduced and has at most simple singularities, that is $W$ is normal with only rational double points as singularities. Hence $W$ is the Galois closure of $\overline \gamma$.
\begin{proof}
Thanks to Proposition \ref{proposition ORDER THREE}, the general fiber $F$ does not contain any inflection point. Moreover, $\ba\cdot F=10$ and by Proposition \ref{proposition CONFIGURATIONS} we know scheme-theoretically the intersection. Hence the divisor $\ba_{|F}$ consists of ten distinct points.
Observe moreover that, as $\ba\cdot E_{k}=0$, it  does not contain any vertical component with respect to $\overline f$.
We can thus conclude that $\ba$ is reduced. This is equivalent to $W$ being a normal surface (see \cite[Proposition 1.1]{Persson81}).

From Proposition \ref{proposition CONFIGURATIONS}, we see that locally $\ba$ has intersection multiplicity at most $2$ with any fiber. This implies that it can have at most double points, i.e. all the possible singularities of $\ba$ are simple points of type $A_n$. These singularities of the branch locus give rise to rational double points of $W$ (see \cite{Persson78}).

As the Galois closure of $\overline \gamma$ is  by definition a normalization of $W$, the last statement is straightforward.
\end{proof}
\end{proposition}

\medskip
\section{The $\lpr$ surfaces and their invariants}\label{section INVARIANTS}

\begin{definition}\label{definition LPR}
Let $(S,[\mathcal L])\in \mathcal W(1,2)$ be such that any element of $|\mathcal L|$ is irreducible. Let $X$ be the minimal model of the Galois closure $W$ of ${\gamma}\colon{S}\dasharrow Y$. We shall call this surface the {\em $\lpr$ surface associated to $(S,\mathcal L)$}.
\end{definition}
This section is devoted to compute the Chern invariants of an $\lpr$ surface $X$ (Theorem \ref{theorem TOPOLOGICAL INDEX} and Theorem \ref{theorem IRREGULARITY}). Furthermore, we shall prove that the Albanese variety of $X$ is the product $S\times S$ and that $X$ is a Lagrangian surface (see Theorem \ref{theorem ALBANESE}).

By abuse of notation, let $E_0,G_{jk},F\subset W$ be the pullbacks of $E_0,G_{jk},F\subset \overline{S}$ via $\alpha_1$. For $1\leq k\leq 3$, the curve $E_k$ does not meet the branch locus $\ba$, hence its pullback consists of two curves $E_k'$ and $E_k''$. By Proposition \ref{proposition Balpha SINGULARITIES}, we are able to compute explicitly the invariants of $X$, as follows.

\begin{theorem}\label{theorem TOPOLOGICAL INDEX}
The $\lpr$ surfaces are of general type, with invariants
$$
K_{X}^2= 198\quad c_2(X)=102\quad \chi(\oo_{X})=25
$$
\begin{proof}
Let us suppose that $\ba$ - and hence $W$ - is smooth. The surface $W$ is minimal. Indeed the $-1$-curves on $W$ come either from the $-1$-curves $L\subset \overline{S}$ such that $\ba\cap L =\emptyset$, or from the $-2$-curves on $\overline{S}$ entirely contained in the branch divisor $\ba$. The only $-1$-curves on $\overline{S}$ are $E_0$ and the $G_{jk}$'s, but they intersect $\ba$. On the other hand, $\overline{S}$ does not contain any $-2$-curve.

Setting $X=W$, the formulas to compute the invariants of $X$ are the following \cite{Persson78}:
$$
K^2_X=2\left(K^2_{\overline{S}}+2p_a(\ba)-2\right)-\frac{3}{2}B^{\cdot 2}_{\alpha_1}
$$
$$
c_2(X)=2c_2(\overline{S})+2p_a(\ba)-2
$$
$$
\chi(\oo_{X})=2\chi(\oo_{\overline{S}})+\frac{p_a(\ba)-1}{2}-\frac{B^{\cdot 2}_{\alpha_1}}{8}\,,
$$
where $p_a(\ba)$ denotes the arithmetic genus of the branch curve. We note that $\overline{S}$ is obtained by blowing up ten times the abelian surface $S$, hence $c_2(\overline{S})=10$ and $\chi(\oo_{\overline{S}})=0$. Moreover, by the adjunction formula and (\ref{equation CLASS OF Balpha}) we have
$$
2p_a(B_{\alpha_1})-2=(K_{\overline S}+B_{\alpha_1}\cdot B_{\alpha_1})=82,
$$
that is $p_a(\ba)=42$. Applying equations (\ref{equation CANONICAL DIVISOR}) and (\ref{equation CLASS OF Balpha}) to the above formulas we compute the invariants. 
By the Enriques-Kodaira classification, $X$ is a surface of general type.

Now, let us assume that $\ba$ is singular. By Proposition \ref{proposition Balpha SINGULARITIES} it has only simple singularities. Consider the canonical resolution of the double covering $\alpha_1\colon W\longrightarrow\overline{S}$ \cite[Section III.7]{BHPV}. We have the following diagram
$$
\xymatrix{ X' \ar[r] \ar[d] & W \ar[d]^{\alpha_1} \\ \overline{S}\,' \ar[r] & \overline{S}, }
$$
where $\overline{S}\,'$ is obtained by blowing up $\overline{S}$ in order to perform the embedded resolution of $B_{\alpha_1}\subset \overline S$ and $X'$ is the smooth surface obtained as a double covering of the strict transform of $\ba$. Following the study in \cite[Section III.7]{BHPV} we see that $X'$ does not contain $-1$-curves, thus $X'$ is the minimal desingularization of $W$. Finally, the formulas to compute the Chern invariants of $X$ are the same we used above \cite[Theorem III.7.2 and Section V.22]{BHPV}.
\end{proof}
\end{theorem}
\begin{remark}\label{family}
We can develop our construction of the Galois closure $W$ from $(S,[\mathcal L])$ in families, thus obtaining a flat three-dimensional family of surfaces of general type with at most canonical singularities over a Zariski open subset of $\mathcal W(1,2)$, contained in the moduli space of surfaces of general type with the above invariants. It would be interesting to understand whether these surfaces form a connected component of this moduli space or not.
\end{remark}

\begin{remark}\label{remark QUOZIENTE A3}
To make our description more complete, let us spend a few words about the surface $\Sigma$ obtained as the quotient of $X$ by the alternating subgroup $A_3\subset S_3$. It can be easily proved that the action of $A_3$ on $X$ has only isolated fixed points. These points correspond to the $80$ points of order three in $S$, and to the points $p\in S$ such that the canonical image of the corresponding curve in $\vert \mathcal L\vert$ has an inflection point of order $4$ in $p$.
The surface  $\Sigma$ is thus a normal surface with rational double points of type $A_2$.  It can be seen as a  double covering $\delta \colon \Sigma\longrightarrow Y$ with branch divisor $B_{\delta}=\overline \gamma_*R_{\overline \gamma}$. The numerical class of this branch is $4C_0+30f$. By arguing as above we can compute the invariants of its canonical resolution $\widetilde \Sigma$, which are the following
$$
K_{\widetilde \Sigma}^2=66, \qquad c_2(\widetilde \Sigma)=258, \qquad \chi(\mathcal O_{\widetilde \Sigma})=27.
$$
Hence $\widetilde \Sigma$ is a surface of general type with an induced stable genus $4$ fibration over $l\colon \widetilde \Sigma\longrightarrow \pr^1$, whose slope $s_l=(K_{\widetilde \Sigma}^2+24)/(\chi(\mathcal O_{\widetilde \Sigma})+27)=3$ reaches the minimum in the slope inequality \cite{C-H}.
Moreover, the irregularity of $\Sigma$ - and hence of $\widetilde \Sigma$ - is 0. Indeed, by the theory of cyclic coverings, $\delta_*\mathcal O_\Sigma=\mathcal O_Y\oplus \mathcal D^{-1}$, where $\mathcal D\cong 2C_0+15f$ is the line bundle associated to the double cover $\delta\colon \Sigma\longrightarrow Y$.
From the Leray spectral sequence we have that
$$
H^1(\oo_\Sigma)=H^1(\delta_*\oo_\Sigma)=H^1(\oo_Y)\oplus H^1(\mathcal D^{-1}).
$$
As $\mathcal D$ is ample, Kodaira vanishing Theorem implies that $H^1(\mathcal D^{-1})=\{0\}$.\\
\end{remark}

We now want to study the space of holomorphic 1-forms $H^0(X, \Omega_X^1)$. It carries a natural $S_3$-representation structure, because the monodromy group $M(\gamma)\cong S_3$ acts on $X$. We recall that the irreducible representations of $S_3$ are the trivial, the anti-invariant, and the standard one \cite[Section 1.3]{FH}, that we denote  $U$, $U'$ and $\Gamma$ respectively.

As in (\ref{diagram XZY}), let us denote by $\pi_i\colon X\longrightarrow \overline{S}$ the composition of the desingularization map $X\longrightarrow W$ with the morphisms $\alpha_i\colon W\longrightarrow \overline{S}$, for $1\leq i\leq 3$. Thus we have the following diagram
\begin{equation}\label{diagram XSY}
\xymatrix{ X\ar[dr]\ar[d]_{\pi_i} & \\ \overline{S}\ar[r]^{\gamma} & Y\\}\,.
\end{equation}
Let $\{\eta_1,\eta_2\}$ be a basis of $H^0(\overline{S},\Omega^1_{\overline{S}})$ and for $1\leq i\leq 3$ and $j=1,2$, let us define the 1-forms
$$
\omega^i_j:=\pi_i^*(\eta_j)\in H^0(X,\Omega^1_X).
$$
Notice that the action of $S_3$ on the $\omega^i_j$'s is given by $\sigma \omega^i_j=\omega^{\sigma(i)}_j$ for any $\sigma\in S_3$. The following result is a consequence of Proposition \ref{proposition p-FORMS}.
\begin{lemma}\label{lemma IRREGULARITY}
The form $\omega\in \bigwedge^2 H^0(X,\Omega^1_X)$ given by
\begin{equation}\label{equation OMEGA}
\omega := \omega^1_1\wedge \omega^1_2 + \omega^2_1\wedge \omega^2_2 + \omega^3_1\wedge \omega^3_2= \frac{3}{2}\,\omega^3_1\wedge \omega^3_2+\frac{1}{2}\,(\omega^1_1-\omega^2_1)\wedge (\omega^1_2-\omega^2_2)
\end{equation}
has rank 4 and it belongs to the kernel of \,$\psi_2 \colon \bigwedge^2 H^0(X,\Omega^1_X)\longrightarrow H^0(X,\Omega^2_X)$. Moreover, the decomposable forms $\omega^3_1\wedge\omega^3_2$ and $(\omega^1_1-\omega^2_1)\wedge (\omega^1_2-\omega^2_2)$ are non-zero in $H^0(X,\Omega^2_X)$.
\begin{proof}
Note that the 1-form $\omega^1_j+\omega^2_j+\omega^3_j$ is $S_3$-invariant for any $j=1,2$. Since the quotient $X/S_3$ is the rational surface $Y=\mathbb{F}_3$ - which does not admit non-trivial 1-forms - we deduce that $\omega^1_j+\omega^2_j+\omega^3_j=0$ on $X$ and hence $\omega^3_j=-\omega^1_j-\omega^2_j$. Thus it is easy to check the equality in (\ref{equation OMEGA}). Then by Proposition \ref{proposition p-FORMS}, we have that $\omega$ has rank 4 and $\omega\in\ker \psi_2$.

The form $\omega^3_1\wedge\omega^3_2$ is the pullback via $\pi_3$ of the $2$-form $\eta_1\wedge\eta_2$ on $\overline{S}$. Since $\eta_1\wedge\eta_2$ is a generator of the one-dimensional vector space $H^0(\overline{S}, \Omega^2_{\overline{S}})$, we have that $\omega^3_1\wedge \omega^3_2$ provides a non-zero holomorphic $2$-form on $X$, that is $\omega^3_1\wedge \omega^3_2\not\in \ker\psi_2$. Therefore $(\omega^1_1-\omega^2_1)\wedge (\omega^1_2-\omega^2_2)\not\in \ker\psi_2$ as well, because $\omega\in\ker\psi_2$.
\end{proof}
\end{lemma}

\begin{theorem}\label{theorem IRREGULARITY}
Let $X$ be a $\lpr$ surface. The space of holomorphic 1-forms on $X$ is $H^0(X, \Omega_X^1)\cong\Gamma\oplus \Gamma$ and hence $q(X)=4$ and $p_g(X)=28$.
\begin{proof}
The vector space $H^0(X, \Omega_X^1)$ admits a decomposition into the direct sum of the irreducible representations of $S_3$ (see \cite[Proposition 1.8 p. 7]{FH}), that is
\begin{equation}\label{equation DECOMPOSIZIONE}
H^0(X, \Omega_X^1)=U^{\oplus a}\oplus U'^{\oplus b} \oplus \Gamma^{\oplus c} \quad \mbox{ for some }\,\,a,b,c\in \mathbb N.
\end{equation}
We want to prove that $a=b=0$ and $c=2$. Since $Y=X/S_3$ is a rational surface, ${h^0(Y,\Omega_Y^1)=0}$ and hence $H^0(X, \Omega_X^1)$ does not contain any invariant element. Thus we have $a=0$.

By Theorem \ref{theorem qCOPIE} the space $H^0(X, \Omega_X^1)$ contains $2$ copies of the standard representation, $\Gamma_{\eta_j}:=\langle \omega^1_j,\omega^2_j,\omega^3_j\rangle$ (following the notation of the Theorem). This implies that $c\geq 2$.
Notice that each $\Gamma_{\eta_j}$ splits with respect to the action of $(12)\in S_3$ in an invariant space $\langle \omega^1_j+\omega^2_j \rangle $ and an anti-invariant one, generated by the one form $\nu_j:=\omega^1_j-\omega^2_j$.  From this we see that  $\dim H^0(X, \Omega_X^1)^{(12)}=c$. From the identification  $H^0(X, \Omega_X^1)^{(12)}=H^0(\overline S, \Omega_{\overline S}^1)=\langle \eta_1,\eta_2\rangle$ we conclude that $c=2$.

Let us prove that $b=0$. Suppose by contradiction that there exists a $1$-form ${\nu_3\in H^0(X,\Omega_X^1)}$ belonging to the anti-invariant representation $U'$ of $S_3$. Let us consider the vector space $R:=\langle \nu_1, \nu_2, \nu_3\rangle$. Since $\nu_1, \nu_2, \nu_3$ are all anti-invariant under the action of $(12)$, we have that $\nu_1\wedge \nu_3$, $\nu_2\wedge \nu_3$ and $\nu_1\wedge \nu_2$ are $\langle (12)\rangle$-invariant. We recall that $X/\langle (12)\rangle=\overline{S}$ and $h^0(\overline{S},\Omega^2_{\overline{S}})=1$. Moreover, $\nu_1\wedge \nu_2\not\in \ker\psi_2$ by Lemma \ref{lemma IRREGULARITY}. Thus the image of the map
$$
\overline{\psi}:=\psi_{2|\wedge^2 R}\colon \wedge^2 R \longrightarrow H^0(X,\Omega_X^2)
$$
is one-dimensional. Therefore $\ker\overline{\psi}$ has dimension 2. We then consider the subspaces ${\langle \nu_1\wedge \nu_2, \nu_1\wedge \nu_3\rangle}$ and $\langle \nu_2\wedge \nu_1, \nu_2\wedge \nu_3\rangle $ of $\bigwedge^2 H^0(X, \Omega^1_X)$. Their intersection with $\ker \overline\psi$ has necessarily dimension one. So, there exist $s,t,w,z\in \mathbb{C}$ such that $\nu_1\wedge (s\nu_3+t\nu_2),\nu_2\wedge (w\nu_3+z\nu_1)\in \ker\,\overline{\psi}$. In particular, there exists a rational function $h$ on $X$ such that $\nu_1=h(s\nu_3+t\nu_2)$ and hence $\nu_2\wedge (w\nu_3+z\nu_1)=\nu_2\wedge (w+h\,zs)\nu_3\in \ker\overline{\psi}$. Then we have that $\nu_2=h_2\nu_3$ for some rational function $h_2$ on $X$. Analogously, there exists $h_1\in K(X)$ such that $\nu_1=h_1\nu_3$. Thus $\nu_1\wedge\nu_2\in\ker \overline{\psi}$, a contradiction. Therefore we have $b=0$.

Thus $V=H^0(X,\Omega^1_X)=\Gamma_{\eta_1}\oplus \Gamma_{\eta_2}$. In particular, we deduce $q(X)=h^0(X,\Omega^1_X)=4$ and $p_g(X)=\chi(\oo_{X})+q(X)-1=28$.
\end{proof}
\end{theorem}

We now study the Albanese variety of $X$. We need to analyze some special fibers of the induced fibration $h:=\overline f\circ \alpha_1\colon W\longrightarrow \pr^1$. Let $D\subset S$ be a hyperelliptic smooth element of $|\mathcal{L}|$. With the same notation of Proposition \ref{proposition BLOW UP}, let $F=D\cup G$ be the corresponding fiber of $\overline  f\colon \overline S \longrightarrow \pr^1$. So, $G$ is a copy of $\pr^1$ attached to $D$ in one node. Let $H\subset W$ denote the pullback of $F\subset \overline{S}$ via $\alpha_1$.
\begin{lemma}\label{lemma TRIANGLE}
The fiber $H$ has three irreducible components $D_1$, $D_2$, $D_3$, that are all copies of $D$ attached two by two in one node as in Figure \ref{figure TRIANGLE} below. Moreover, for a suitable choice of the indices, ${\alpha_i}_{|D_i}\colon D_i\longrightarrow G$ is the hyperelliptic map, whereas for $i\neq j$, ${\alpha_j}_{|D_i}\colon D_i \longrightarrow D$ is either the identity map or the hyperelliptic involution.
\end{lemma}
\begin{proof}
The first part of the statement follows from the description of $B_{\alpha_1|F}$ given in Proposition \ref{proposition CONFIGURATIONS}. Indeed, $B_{\alpha_1|G}$ consist of the $8$ Weierstrass points of the hyperelliptic map ${D\longrightarrow G\cong\pr^1}$, whereas $B_{\alpha_1|D}$ is the point $p_0$ with multiplicity $2$. Thus the inverse image $D_1$ of $G$ is a copy of $D$, whereas the inverse image of $D$ is given by two copies of $D$ attached in one node.

The second statement follows from the definition of $W\subset \overline S\times \overline S$ and of the $\alpha_i$'s. Indeed, observe that, given a general point $q\in G\cong\pr^1$, its preimages via $\overline \gamma$ are $q\in G$ itself, and $q_1, q_2\in D$ (the two preimages of the hyperelliptic involution). The corresponding fibers of $\alpha_1$ are $a=(q,q_1)$ and $b=(q,q_2)$. So, $\alpha_2(a)=q_1$, $\alpha_2(b)=q_2$, and $\alpha_3(a)=q_2$, $\alpha_3(b)=q_1$, as claimed.
\end{proof}
\begin{center}
\begin{figure}[h!]
\input FIBREtriangle.pstex_t
\caption{}
\label{figure TRIANGLE}
\end{figure}
\end{center}

\begin{theorem}\label{theorem ALBANESE}
Let $X$ be a $\lpr$ surface associated to $(S,\mathcal L)$. The Albanese variety of $X$ is $ \mathrm{Alb} (X)=S\times S$. Furthermore, $X$ is a Lagrangian surface.
\begin{proof}
The universal property of the Albanese morphism induces a morphism of abelian varieties $\theta\colon \Alb (X)\longrightarrow S\times S$ which is an isogeny because $\dim \Alb (X)=q(X)=4$ by Theorem \ref{theorem IRREGULARITY}. Thus we have an induced inclusion in singular homology
$$
\theta_*\colon H_1(X,\mathbb{Z})=H_1(\Alb (X),\mathbb{Z})\longrightarrow H_1(S\times S,\mathbb{Z})=H_1(S,\mathbb{Z})\times H_1(S,\mathbb{Z}).
$$
We want to prove that $\theta_*$ is surjective, i.e. that $\theta$ is an isomorphism of abelian varieties.
Let us consider one of the fibers of $h:= \overline{f}\circ \alpha_1$ inspected in Lemma \ref{lemma TRIANGLE}, ${H=D_1\cup D_2 \cup D_3\subset W}$, and let $\widetilde{H}\subset X$ be its pullback via the normalization morphism. By Proposition \ref{proposition Balpha SINGULARITIES}, the only possible singularities of $W$ are rational double points. Hence $\widetilde{H}$ may differ from $H$ only for some $-2$-curves, that are contracted by the Albanese morphism $a\colon X\longrightarrow \Alb(X)$. So, let ${\widetilde{H}=D_1\cup D_2 \cup D_3\cup R}$, where the $D_i$'s are the strict transforms on $X$ of the components of $H$, and $R$ is the - possibly null - divisor given by the $-2$-curves.

Thanks to Lemma \ref{lemma TRIANGLE}, we have $\alpha_1(D_2)=D\subset \overline S$, whereas $\alpha_2(D_2)=G \cong \pr^1\subset \overline S$. Thus the image of $D_2$ in $S\times S$ is $D\times \{0\}$. By Proposition \ref{proposition SMOOTH}, for any smooth member $C\in |\mathcal{L}|$, $S$ is naturally identified with $J(C)/\pi^*E$, where $\pi\colon C\longrightarrow E$ is the bielliptic involution. So, $S$ fits in the following diagram of abelian varieties
\begin{equation}\label{BAA}
1\longrightarrow E\longrightarrow J(D)\stackrel{\zeta}{\longrightarrow} S\longrightarrow 1.
\end{equation}
The image of the composition
$$
H_1(D_2,\mathbb{Z})\longrightarrow H_1(X,\mathbb{Z})\stackrel{{\theta}_*}{\longrightarrow}H_1(S,\mathbb{Z})\times H_1(S,\mathbb{Z})
$$
is $H_1(S,\mathbb{Z})\times \{0\}$. Indeed, it can be naturally identified with the homomorphism induced in homology by the sequence (\ref{BAA}):
$$
H_1(D,\mathbb{Z})=H_1(J(D),\mathbb{Z})\longrightarrow H_1(S,\mathbb{Z}),
$$
which is surjective because the map $\zeta$ in (\ref{BAA}) has connected fibers. The same argument applied to $D_1$ proves that the image of ${\theta}_*$ contains $\{0\} \times H_1(S,\mathbb{Z})$ as well. Thus $\theta_*$ is surjective as wanted.

In particular, $X$ admits a generically finite morphism of degree one into the four-dimensional abelian variety $S\times S$ and by Lemma \ref{lemma IRREGULARITY} there exists a holomorphic $2$-form ${\omega\in \bigwedge^2H^0(X,\Omega^1_X)=H^{2,0}(\Alb (X))}$ of rank $4$ such that $\omega\in \ker \psi_2$. Thus $X$ is a Lagrangian surface.
\end{proof}
\end{theorem}

\medskip
\section{The Galois closure of $\overline{\gamma}$ is non-fibred}\label{section NON-FIBRED}

We are going to prove that for a general choice of the abelian surface $S$ with the $(1,2)$-polarization, the associated $\lpr$ surface $X$ does not admit fibrations on curves of genus $\geq 2$.

\begin{proposition}\label{non fibrata >2}
Let $(S,\mathcal [L])\in \mathcal W(1,2)$ be such that $\vert \mathcal L\vert$ contains only irreducible elements.
Then the associated $\lpr$ surface does not admit fibrations on curves of genus $\geq 3$.
\begin{proof}
Let $X$ be the $\lpr$ surface associated to  $(S,\mathcal [L])$. Suppose by contradiction that $\mu \colon X\longrightarrow T$ is a fibration over a smooth curve of genus $g\geq 3$. Clearly $g$ is at most equal to the irregularity $q(X)=4$. As $X$ has finite Albanese morphism, it cannot be $g=4$.

Suppose that $g=3$. Let $\tau_1\colon X\longrightarrow X$ the involution associated to the double covering $\alpha_1\colon X\longrightarrow \overline S$, and  consider the fibration $\nu=\mu\circ \tau_1\colon X\longrightarrow T$. Consider the images of the induced  pushforward maps
$$
W= \mu^*( H^{1,0}(T))\subseteq H^{1,0}(X),\quad
V=\nu^*(H^{1,0}(T))\subseteq H^{1,0}(X)
$$
The spaces $V$ and $W$ are $3$-dimensional and do not coincide, because $\nu$ and $\mu$ are different fibrations.
Hence, the intersection $V\cap W$ has dimension $2$. Let $u_1, u_2$ be a basis of $V\cap W$. As $u_i\in V$, we have that $u_i\wedge v=0$ for any $v\in V$, and similarly $u_i\wedge w=0$ for any $w\in W$. But $V+W=H^{1,0}(X)$, and so  that the cup homomorphism $\psi_2\colon \bigwedge^2 H^{1,0}(X) \longrightarrow H^{2,0}(X)$ should be trivial, which is absurd.
\end{proof}
\end{proposition}

It remains to exclude the case of a fibration over a genus $2$ curve. The argument we develop for this case is less straightforward and it is based on the properties of the six hyperelliptic elements of the linear pencil $\vert \mathcal L\vert$.
We will need the following simple result.
Let $\mathcal M_3$ be the moduli space of smooth genus $3$ curves; and let $\mathcal Y_3$ and $\mathcal B_3$ be the hyperelliptic and the bielliptic locus respectively.
\begin{lemma}\label{moduli}
If $(S, [ \mathcal L])$ is general in $\mathcal W(1,2)$, then the hyperelliptic elements in $\vert \mathcal L\vert$ are general in $\mathcal Y_3\cap \mathcal B_3$.
\begin{proof}
Let $[C]\in \mathcal Y_3\cap \mathcal B_3$ be a curve with a unique bielliptic involution, and consider the associated double covering $\pi\colon C\longrightarrow E$. Recall that the generalized Prym variety $\prym =J(C)/\pi^*E$ is an abelian surface with a $(1,2)$-polarization.
By associating $\prym$ to $C$, we construct a map
$$
\Phi \colon \mathcal Y_3\cap \mathcal B_3 \,-\!\!\dasharrow \mathcal W(1,2).
$$
For any general $(S, [ \mathcal L])\in\mathcal W(1,2)$, the linear series $|\mathcal{L}|$ possesses exactly six smooth hyperelliptic curves. Moreover, $\dim \mathcal Y_3\cap \mathcal B_3=3=\dim \mathcal W(1,2)$, hence $\Phi $ is generically finite dominant map and the assertion follows.
\end{proof}
\end{lemma}

\begin{proposition}\label{non fibrata 2}
For a general choice of $(S,[\mathcal L])\in \mathcal W(1,2)$, the associated $\lpr$ surface is not fibred over curves of genus $2$.
\begin{proof}
Let us suppose by contradiction that $\mu \colon X\longrightarrow T$ is a fibration over a smooth curve of genus $2$.
Let us consider again  the ``triangular'' fibers of Lemma \ref{lemma TRIANGLE}.
Recall that for any hyperelliptic element $D\in |\mathcal L|$ there is one of these fibers $H=D_1\cup D_2\cup D_3\subset X$, whose irreducible components are all smooth hyperelliptic and bielliptic curves of genus $3$ isomorphic to $D$.
At least one of the components is not contracted by $\mu$, otherwise the map $\mu$ should factor through the fibration $h\colon X\longrightarrow \pr^1$, which is of course impossible.
So, there is a finite morphism $D\longrightarrow T$. By the Hurwitz formula, this morphism has to be \'etale of degree $2$.
Consider the hyperelliptic involution $\iota$ and the bielliptic one $j$ over $D$. The composition $\sigma =\iota\circ j$ is a fixed-point-free involution, and - by Hurwitz formula again - the quotient $D/\langle \sigma\rangle$ is smooth of genus two.

A general curve  in $\mathcal Y_3\cap \mathcal B_3$ has automorphism group isomorphic to  $ \Z/2\Z\times \Z/2\Z$, generated by the hyperelliptic and the bielliptic involutions.
By Lemma \ref{moduli}, under our assumptions we can suppose that $D$ is general in $\mathcal Y_3\cap \mathcal B_3$, and hence that  the curve $T$ coincides with $D/\langle\sigma\rangle$.
So, for a general choice of $(S, [\mathcal L])$, the morphism $\mu\colon D\longrightarrow T$ coincides with the quotient map $D\longrightarrow D/\langle\sigma\rangle$.
 Let us consider the induced map on the Jacobians
$$
\vartheta \colon J(D)\longrightarrow J(T).
$$
Let $E= D/\langle j\rangle$. We shall now prove that
$\pi^*E\subseteq \ker \vartheta$.
Recall that the hyperelliptic involution $\iota$ induces a permutation of $\{x_0,\ldots ,x_3\}$, which are the fixed points of the bielliptic involution $j$ (cf. Proposition \ref{proposition HYPERELLIPTIC FIBERS}). Let $e$ be an element of $\pi^*E\subset J(D)$. Then by Proposition \ref{proposition SMOOTH}
$$e\sim_D s +j(s)-2x_0$$
for some $s\in D$. Observe that $s+\iota (s)\sim_Dx_0+\iota(x_0)$, so we have
\begin{equation}\label{gia}
\begin{array}{ll}
e & \sim_Ds +j(s)-2x_0-s-\iota (s) +x_0+\iota(x_0)\\
\, & \sim_D j(s)-\iota(s)+ \iota(x_0)-x_0\\
\, & \sim_D \sigma(\iota (s))- \iota (s) + \sigma(x_0)-x_0.
\end{array}
\end{equation}
Hence, $\pi^*E\subseteq \ker \vartheta$, as wanted. Moreover, note that $\pi^*E= (\ker \vartheta)^\circ$, because of the connectedness of $\pi^*E$, and --from the above equation-- we have that $\pi^*E$ has index two in $\ker \vartheta$.
We thus obtain the following diagram
$$
\xymatrix{
& & 0  \ar[d]  & & \\
& & \ker \vartheta \ar[d] & &\\
0 \ar[r] &\pi^*E\ar[ru]\ar[r] & J(D) \ar[r]  \ar_{\vartheta}[d]  &S\ar^{\theta}[dl] \ar[r] & 0\\
& & J(T)\ar[d] & &\\
& & 0 & &
}
$$
where $\theta\colon S\longrightarrow J(T)$ is a degree $2$ isogeny of abelian variety with kernel $\Z/2\Z\cong \langle x_0-\sigma(x_0) \rangle\subset S$.

For any hyperelliptic element $D\in |\mathcal L|$, we can carry on the above construction.
Recall that for any $i=1,2,3$, there exists an hyperelliptic element in $|\mathcal L |$ whose hyperelliptic involution maps $x_0$ to $x_i$. Hence, for any $i$, we have an isomorphism of principally polarized abelian varieties from $S/\langle x_i-x_0 \rangle$ to $J(T)$.
In the next lemma, we prove that this is impossible. We have thus reached a contradiction and the proof is concluded.
\end{proof}
\end{proposition}
\begin{lemma}
For a general $(S,[\mathcal L])\in \mathcal W(1,2)$ and for $i\neq j$,
the principally polarized abelian surfaces $S/\langle x_i-x_0 \rangle$ and $S/\langle x_j-x_0\rangle $ are not isomorphic.
\begin{proof}
Let $\mathcal W(1,1)$ be the moduli space of principally polarized abelian surfaces. Given $(S,[\mathcal L])\in \mathcal W(1,2)$, by (\ref{tilambda}) we have that $T(\mathcal L)=\{x_0,\ldots,x_3\}$. Hence there exists a well defined morphism
$$
\Psi\colon \mathcal W(1,2)\longrightarrow \sym^3\mathcal W(1,1)
$$
defined as $\Psi((S,[\mathcal L]))=(S/\langle x_1-x_0\rangle, [\mathcal L_1])+(S/\langle x_2-x_0\rangle, [\mathcal L_2])+(S/\langle x_3-x_0\rangle,  [\mathcal L_3])$, where $[\mathcal L_i]$ is the induced polarization. Therefore it suffices to show that there exist a couple ${(S_0, [\mathcal L_0])\in \mathcal W(1,2)}$ whose image under $\Psi$ does not belong to the diagonal.

Let us construct such a polarized surface. Let $E, F$ be two elliptic curves, and consider the abelian surface $S=E\times F$ together with the product polarization given by the line bundle $\mathcal L=\mathcal O_S(2E+F)$. Note that these are the very same surfaces of $\mathcal W(1,2)$ that we excluded in our construction at the beginning of Section \ref{section BARTH}; however, as the map $\Psi$ is everywhere defined over $\mathcal W(1,2)$, there is no contradiction in considering these surfaces in this context.
In this case $T(\mathcal L)$ is $\{x_0=(e_0, 0), \ldots , x_3=(e_3,0)\}\subset E\times F$, where $e_0,\ldots,e_3$ are the four points of order two on $E$. Hence $S/ \langle x_i-x_0\rangle\cong E/\langle e_i-e_0\rangle\times F$. It is now sufficient to choose $E$ such that its quotients by the 2-torsion points are not isomorphic. Let us consider for instance the elliptic curve $E=\bC/ \Gamma$, where $\Gamma =2i\Z\oplus \Z$. Then, setting $e_1=i$, $e_2=1/2$, $e_3=1/2+i$, we have
$$
\frac{E}{\langle e_1-e_0\rangle}=\frac{\bC}{i\Z\oplus \Z},
\quad \frac{E}{\langle e_2-e_0\rangle}=\frac{\bC}{2i\Z\oplus1/2 \Z}\cong \frac{\bC}{4i\Z\oplus \Z},
$$
$$
\frac{E}{\langle e_3-e_0\rangle}=\frac{\bC}{(1/2+i)\Z\oplus 1/2\Z}\cong \frac{\bC}{(1+2i)\Z\oplus \Z}.
$$
It is immediate to check that there is no integer matrix
$\begin{pmatrix}a&b\\
c&d\end{pmatrix}\in SL_2(\Z)$
such that $\frac{ai+b}{ci+d}=4i$, hence the first two curves are not isomorphic.
\end{proof}
\end{lemma}

The existence of a fibration with base curve of genus $\geq 2$ on a surface is a topological condition, as proved by Siu \cite{Siu}, Beauville and Catanese \cite{Catanese}.
Thus, from Proposition \ref{non fibrata >2}  and \ref{non fibrata 2} we can derive the following conclusion.
\begin{theorem}\label{nf}
For any choice of $(S,[\mathcal L])\in \mathcal W(1,2)$, the associated $\lpr$ surface is not fibred over a curve of genus $\geq 2$.
\end{theorem}

As mentioned in the introduction, the knowledge of the cup product homomorphism
$$\rho_2\colon H^1(X,\C)\wedge H^1(X,\C) \longrightarrow H^2(X,\C)$$
for a $\lpr$ surface $X$
would allow to compute completely the nilpotent tower of $\pi_1(X)$. As a first step we compute the dimension and the decomposition as $S_3$-representation of $\ker \rho_2$.

\begin{proposition}\label{kerrho2}
For any $\lpr$ surface $X$ we have that $\ker \rho_2=\Gamma\oplus U^{\oplus 5}$ as an $S_3$-representation.
\begin{proof}
Recall the Hodge decomposition of $\rho_2= \rho_2^{2,0}\oplus  \rho_2^{1,1}\oplus {\rho_2^{0,2}}$, where $\rho_2^{2,0}=\psi_2$, and $\rho_2^{0,2}$ is its conjugate.
By Theorem \ref{nf} and the Castelnuovo-de Franchis Theorem, we have that $\ker \psi_2=U$.
Le us study $\ker \rho_2^{1,1}\colon H^{1,0}(X)\otimes H^{0,1}(X)\longrightarrow H^{1,1}(X)$.
From Theorem \ref{theorem IRREGULARITY}, we deduce that, as a representation,
$$H^{1,0}(X)\otimes H^{0,1}(X)=\Gamma^{\oplus 4}\oplus U'^{\oplus 4}\oplus U^{\oplus 4}.$$
Let  $g\colon X\longrightarrow \mathbb F_3$ be the quotient map via the action of $S_3$. The space $g^\ast H^{1,1}(\mathbb F_3)\subset H^{1,1}(X)$ has dimension $2$ and it is generated by the class of a fiber $[g^\ast f]$  and by the class of the negative section $[g^\ast C_0]$.
Note that $\rho_2^{1,1}$ coincides with the pushforward map induced by the Albanese morphism $a\colon X\longrightarrow \Alb (X)$
$$a^\ast \colon H^{1,1}(\Alb(X))\longrightarrow H^{1,1}(X).$$
Note that  $g^\ast C_0$ consists of $-3$-curves, so it is contracted by $a$. Hence  we have that $\mbox{Im}\rho_2^{(1,1)}\cap g^\ast H^{1,1}(Y)$ is $1$-dimensional.
We can therefore conclude that there is an invariant $3$-dimensional space $U^{\oplus 3}$ contained in the kernel.

Consider now the part of $\rho_2^{(1,1)}$ invariant by a transposition, e.g. $(12)\in S_3$
$$U^{\oplus 4}\oplus V=\left( H^{1,0}(X)\otimes H^{0,1}(X)\right)^{(12)}\longrightarrow \left(H^{1,1}(X)\right)^{(12)}= H^{1,1}(\overline S),$$
where $V$ is a $4$-dimensional $(12)$-invariant space coming from the copies of the standard representation $\Gamma^{\oplus 4}\subset H^{1,0}(X)\otimes H^{0,1}(X)$.
The image of this morphism has dimension $h^{1,1}(S)=4$; as only $3$ copies of $U$ are contained in the kernel, it follows that $\ker \rho_2^{1,1}$ has non-trivial intersection with $V$, and we conclude that there is at least one $\Gamma$ contained in it. In \cite[Theorem 1]{PC}, it is proved that $\dim \ker \rho_2\leq 7$, hence the statement follows.
\end{proof}
\end{proposition}

Note that  the above result  proves that the estimate in \cite{PC} is sharp.

\medskip
\section{Monodromy of 3-torsion points of $S$ and smoothness of the Galois closure}\label{section MONODROMY}

In this section we are going to prove that the Galois closure $W\subset \overline{S}\times \overline{S}$ constructed over a general point $(S,[\mathcal{L}])\in \mathcal{W}(1,2)$ is smooth. This fact will be important to determine the branch locus of the Albanese map for the $\lpr$ surfaces.

To start, we recall some basic facts about moduli of $(1,2)$-polarized abelian surfaces following \cite[Chapter 8]{L-B}. Let $\mathcal H_2$ denote the \emph{Siegel upper half space}, which provides a moduli space for polarized abelian surfaces of type $(1,2)$, together with a symplectic basis. Let ${D=\begin{pmatrix}1&0\\ 0&2\end{pmatrix}}$ be the matrix associated to the polarization and let us consider the symplectic group of matrices preserving the intersection form, that is
$$
\Gamma_D:=Sp_4^{(1,2)}(\Z)=\left \{ R\in M_4(\mathbb{Z}) \mid  R\begin{pmatrix}0&D\\-D&0\end{pmatrix} R^t=\begin{pmatrix}0&D\\-D&0\end{pmatrix}\right\}.
$$
The group $\Gamma_D$ acts on $\mathcal H_2$ and the moduli space $\mathcal W(1,2)$ is the quotient $\mathcal H_2/\Gamma_D$ parametrizing the isomorphism classes of $(1,2)$-polarized abelian surfaces. In particular, given a $(1,2)$-polarized abelian surface $A=\mathbb{C}^2/\Lambda$ with a basis $\left\{\lambda_1,\lambda_2,\mu_1,\mu_2\right\}$ of the lattice $\Lambda$ for $D$, any element of $\Gamma_D$ can correspond to an isomorphism $\varphi\colon \left(A; \lambda_1,\lambda_2,\mu_1,\mu_2 \right)\longrightarrow \left(A'; \lambda'_1,\lambda'_2,\mu'_1,\mu'_2 \right)$, consisting of a change of the basis on $\Lambda$.

We recall that for any point $\left(A; \lambda_1,\lambda_2,\mu_1,\mu_2 \right)\in \mathcal H_2$, the set $\left\{\overline{\frac{1}{3}\lambda_1},\overline{\frac{1}{3}\lambda_2},\overline{\frac{1}{3}\mu_1},\overline{\frac{1}{3}\mu_2}\right\}$ is a basis of the group of $3$-torsion points of $A$, where the bar means that we are taking the equivalence classes modulo $\Lambda$. Let $\Gamma'_D\subset \Gamma_D$ be the subgroup of the matrices representing the isomorphisms of $(1,2)$-polarized abelian surfaces $\varphi\colon \left(A; \lambda_1,\lambda_2,\mu_1,\mu_2 \right)\longrightarrow \left(A'; \lambda'_1,\lambda'_2,\mu'_1,\mu'_2 \right)$ such that $\varphi(\overline{\frac{1}{3}\lambda_1})=\overline{\frac{1}{3}\lambda'_1}$. Therefore we can define the quotient $\mathcal{F}:=\mathcal H_2/\Gamma'_D$ and we have the following diagram
$$
\xymatrix{
\mathcal H_2  \ar_{\delta}[dd]\ar^{\delta_1}[rd]  &  \\
 & \mathcal{F}\ar^{\delta_2}[ld]\\
\mathcal{W}(1,2) &
}
$$
where the map $\delta_2$ is finite (cf. \cite[Theorem 8.3.1]{L-B}). We would like to note that the points of $\mathcal{F}$ are isomorphism classes of $(1,2)$-polarized abelian surfaces together with a torsion point of order $3$, that is the point of coordinates $\left(\frac{1}{3},0,0,0\right)$ with respect to the basis $\left\{ \lambda_1,\lambda_2,\mu_1,\mu_2 \right\}$ of the lattice. The following result is probably well known, but we were not able to find adequate references, hence we give a sketch of the proof.
\begin{theorem}\label{theorem MONODROMY}
Let $(S,[\mathcal{L}])\in \mathcal{W}(1,2)$ be a general point. Then there is a one-to-one correspondence between the set of the $3$-torsion points on $S$ and the fiber of $\delta_2\colon \mathcal{F}\longrightarrow \mathcal{W}(1,2)$ over $(S,[\mathcal{L}])\in \mathcal{W}(1,2)$. Moreover, the monodromy group of $\delta_2$ acts transitively on such a fiber.
\begin{proof}
Let $\left(A; \lambda_1,\lambda_2,\mu_1,\mu_2 \right)\in \mathcal H_2$ be such that
$\delta \left(A; \lambda_1,\lambda_2,\mu_1,\mu_2 \right)=(S,\mathcal{L})$ and let $\delta_1\left(A; \lambda_1,\lambda_2,\mu_1,\mu_2 \right)=\left[A; \overline{\frac{1}{3}\lambda_1} \right]\in \mathcal{F}$. By the construction of $\mathcal{F}$, the degree of $\delta_2$ is at most $80$. Moreover, as $\mathcal{H}_2$ is connected, the quotient $\mathcal{F}$ is connected as well. Therefore the action of the monodromy group of $\delta_2$ is transitive. To conclude the proof, we have then to show that $(S,[\mathcal{L}])$ admits exactly $80$ preimages on $\mathcal{F}$, each corresponding to a different $3$-torsion point. Hence we have to prove that for any $3$-torsion point $p\in A$, there exists an isomorphism of $(1,2)$-polarized abelian surface $\varphi\colon \left(A; \lambda_1,\lambda_2,\mu_1,\mu_2 \right)\longrightarrow \left(A'; \lambda'_1,\lambda'_2,\mu'_1,\mu'_2 \right)$ such that $\varphi(p)=\overline{\frac{1}{3}\lambda'_1}$. By fixing a symplectic basis $\left\{\lambda_1,\lambda_2,\mu_1,\mu_2 \right\}$ for the lattice $\Lambda$, we can identify the group of 3-torsion points on $A$ with the group $\frac{1}{3}\left(\mathbb{Z}/3\mathbb{Z}\right)^4$. Then notice that proving the thesis is equivalent to see that for any non-zero 4-tuple $\left(a_1,a_2,b_1,b_2\right)\in \frac{1}{3}\left(\mathbb{Z}/3\mathbb{Z}\right)^4$, there exists a transformation $\tau$ of the symplectic basis $\left\{\lambda_1,\lambda_2,\mu_1,\mu_2\right\}$ preserving the intersection form such that $\tau\left(\frac{1}{3},0,0,0\right)=\left(a_1,a_2,b_1,b_2\right)$. Then it is easy to see that the group of transformation generated by the following changes of basis
\begin{displaymath}
\begin{array}{lcl} & \tau_1 & \\ \lambda_1 & \mapsto & -\lambda_1\\ \lambda_2 & \mapsto & -\lambda_2\\ \mu_1 & \mapsto & -\mu_1\\ \mu_2 & \mapsto & -\mu_2
\end{array}\quad
\begin{array}{lcl} & \tau_2 & \\ \lambda_1 & \mapsto & \lambda_2\\ \lambda_2 & \mapsto & -\lambda_1\\ \mu_1 & \mapsto & \mu_1\\ \mu_2 & \mapsto & \mu_2
\end{array}\quad
\begin{array}{lcl} & \tau_3 & \\ \lambda_1 & \mapsto & \lambda_1\\ \lambda_2 & \mapsto & \lambda_2\\ \mu_1 & \mapsto & \mu_2\\ \mu_2 & \mapsto & -\mu_1
\end{array}
\end{displaymath}
\begin{displaymath}
\begin{array}{lcl} & \tau_4 & \\ \lambda_1 & \mapsto & \lambda_1+\lambda_2\\ \lambda_2 & \mapsto & \lambda_2\\ \mu_1 & \mapsto & \mu_1\\ \mu_2 & \mapsto & \mu_2
\end{array}\quad
\begin{array}{lcl} & \tau_5 & \\ \lambda_1 & \mapsto & \lambda_1\\ \lambda_2 & \mapsto & \lambda_2\\ \mu_1 & \mapsto & \mu_1+\mu_2\\ \mu_2 & \mapsto & \mu_2
\end{array}
\begin{array}{lcl} & \tau_6 & \\ \lambda_1 & \mapsto & 3\lambda_1-\mu_1\\ \lambda_2 & \mapsto & \lambda_2+\mu_2\\ \mu_1 & \mapsto & \mu_1-2\lambda_1\\ \mu_2 & \mapsto & 3\mu_2+2\lambda_2
\end{array}
\end{displaymath}
acts on $\left(\frac{1}{3},0,0,0\right)$ as wanted, and the proof ends.
\end{proof}
\end{theorem}

\begin{remark}
Note that the above result is essentially due to the fact that the quotient by $\mathbb{Z}/3\mathbb{Z}$ of the $(1,2)$-polarization coincides - up to sign - with the quotient of the principal polarization $(1,1)$. If we  consider torsion points of another order, the situation may be completely different. For instance it is easy to see that the action of $\Gamma_D$ on the 15 points of order 2 separates the three order $2$ base points of $\vert\mathcal{L}\vert$ from the other 12.
\end{remark}

\medskip
We are now ready to prove the following

\begin{theorem}
For a general $(S, [\mathcal L])\in \mathcal{W}(1,2)$, the branch locus $\ba$ is smooth. Therefore the Galois closure $W$ of $\overline{\gamma}$ is smooth and it is the $\lpr$ surface itself.
\begin{proof}
We recall that $B_{\alpha_1}$ is reduced with at most double points as singularities by Proposition \ref{proposition Balpha SINGULARITIES}.
From the study of the intersections of $B_{\alpha_1}$ with the fibers we made in Proposition \ref{proposition CONFIGURATIONS}, we deduce that the possible singular points of $B_{\alpha_1}$ may only lie on $B_{\alpha_1}\cap E_0$ and on the set
$$
\mathcal{B}=\left\{p\in\overline S\setminus E_0 \,| \textrm{ the fiber $F$ through $p$ is smooth non-hyperelliptic and }\,3p+p_0\sim_F K_F \right\}.
$$
A point $p\in \mathcal{B}$ is an order 3 inflection point for the canonical image of the fiber $F$ passing through $p$, with tangent line $\overline{p_0p}$. Note that by Proposition \ref{proposition ORDER THREE}, $\mathcal{B}$ consists of the inverse images of the points of order $3$ in $S$.
Hence the cardinality of $\mathcal{B}$ is $80$. From the monodromy result of Theorem \ref{theorem MONODROMY}, we know that for a general choice of $(S, [\mathcal L])\in \mathcal W(1,2)$ the curve $B_{\alpha_1}$ is either smooth or singular at every point of $\mathcal{B}$. The arithmetic genus of $B_{\alpha_1}$ is $p_a(B_{\alpha_1})=42$ by adjunction formula. We note that the restriction of $\overline f$ to $B_{\alpha_1}$ provides a $10:1$ finite morphism $B_{\alpha_1}\longrightarrow \pr^1$.
As $B_{\alpha_1}$ does not have vertical components, it follows that its irreducible components are at most ten.
If $B_{\alpha_1}$ were singular at any point of $\mathcal{B}$, as any double point increases the arithmetic genus by at least one, the geometric genus would be at most $p_a(B_{\alpha_1})-80\leq -38$, which gives a contradiction. Hence $\ba$ is smooth at $\mathcal{B}$.

Let us now consider the set $E_0\cap B_{\alpha_1}$. It can be naturally divided in the following three subsets.
$$
C_1=\{ p\in E_0\cap B_{\alpha_1} \mid \mbox{ the fiber } F \mbox{ throught } p \mbox{ is singular and irreducible}\}
$$
$$
C_2=\{  p\in E_0\cap B_{\alpha_1} \mid \mbox{ the fiber } F \mbox{ throught } p \mbox{ is reducible}\}
$$
$$
C_3=(E_0\cap B_{\alpha_1})\setminus (C_1\cup C_2).
$$
Note that the cardinality of $C_1$ is $\sharp C_1=12$, whereas $\sharp C_2=6$. By Proposition \ref{proposition Balpha}, we obtain that $B_{\alpha_1}\cdot E_0=22$, therefore
$$
22= B_{\alpha_1}\cdot E_0\geq \sharp C_1+\sharp C_2+ \sharp C_3\geq 18+ \sharp C_3
$$
Hence, for any $(S,[\mathcal L])$, the cardinality of $C_3$ is at most 4. We recall that the points in this set correspond to the inflection points of order 4 lying on the canonical plane image of the fibers. Notice that if $C_3$ consists of exactly $4$ points, then $B_{\alpha_1}$ turns out to be smooth on $E_0$, because in this case it necessarily has local intersection multiplicity $1$ with $E_0$ on $C_1$ and $C_2$. The following provides an example of this situation.
\begin{example}\label{example FERMAT}
Let us consider the plane Fermat curve $C\subset \pr^2$ given by the equation $C\colon x^4+y^4-z^4=0$. It is a smooth curve of genus $3$. The involution $(x:y:z)\mapsto (-x:y:z)$ of $\pr^2$ restricts to a bielliptic involution on $C$, with four fixed points: $c_0=(0: 1:1 )$, $c_1=(0:-1 :1 )$, $c_2=(0:i :1 )$ and $c_3=(0:-i : 1)$. Of course the quotient of $C$ by this involution is a plane elliptic curve $E$, with affine equation $x^2+y^4-1=0$.

Let $\pi\colon C\longrightarrow E$ be the bielliptic quotient map. It is immediate to check that $C$ has inflection points of order $4$ at all the $c_i$'s. Let us consider the abelian surface  $S= \pic^0(C)/\pi^*E$ associated to $C$, where we fix the origin at the image of $c_0$. Then the $(1,2)$-polarization is $\mathcal L\cong \mathcal O_S(C)$. Call $C_i:=C+(c_i-c_0) \subset S$ the translation of $C$ by $c_i$. As $\{c_0,\ldots, c_3\}=T(\mathcal L)$ by the very definition of $T(\mathcal L)$ (cf. (\ref{tilambda})), we have $\mathcal O(C_i)\cong \mathcal L$. So, the $C_i$'s are four distinct elements of $|\mathcal L|$, each having an inflection point of order $4$ at $c_0$.
\end{example}

Thus we exhibit a surface $(S, \mathcal L)$ such that $B_{\alpha_1}$ is smooth on $E_0$. Hence there exists a non-empty Zariski open subset of $\mathcal W(1,2)$ such that for any $(S, [\mathcal L])$ lying on it, the branch curve $B_{\alpha_1}$ is smooth on $E_0$, and the proof is concluded.
\end{proof}
\end{theorem}

\medskip
\section{On a conjecture on the topological index}\label{topological}

Suppose that $\ker \psi_2$ is non-trivial and let $w\in \ker \psi_2$ be a non-zero element.
Let $V\subset H^0(X, \Omega_X^1)$ be the subspace of minimal dimension such that $w\in \bigwedge^2 V$.
Following \cite{BNP} we say that  $X$ is \emph{generalized Lagrangian} if there exists $w\in \ker\psi_2$ of rank $2n$ such that $V$ generically generates
$\Omega_X^1$.
Clearly, a Lagrangian variety is in particular generalized Lagrangian.


Consider the homomorphism $\psi_2$ restricted to $\bigwedge^2 V$, and let $\overline V= \psi_2(\bigwedge^2 V)\in H^0(X, \omega_X)$. Definite $F_V$ as the divisorial part of the base scheme of the linear system $\overline V$. In other words, consider the evaluation morphism $\overline V\otimes \mathcal O_X\longrightarrow \omega_X$ and tensor it with $\omega_X^{-1}$
$$
\overline V\otimes \omega_X^{-1}\longrightarrow \mathcal I\subseteq  \mathcal O_X.
$$
The ideal defining $F_V$ is the double dual $\mathcal I^{**}$. Let $a\colon X\longrightarrow \Alb (X)$ be the Albanese map of $X$. Observe that  any reduced divisor contracted by $a$ is necessarily contained in the base locus of $F_V$, but not vice versa.
We recall from \cite{BNP} that $F_V$ is said to be \emph{contracted by} $V$ if it is reduced and the restriction of the evaluation morphism
$$
V\otimes \oo_X\longrightarrow \Omega_X^1\longrightarrow \Omega_{X|F_V}^1=\omega_{F_V},
$$
vanishes. One of the main results in \cite{BNP} is the following.
\begin{theorem}\cite[Theorem 1.2]{BNP}\label{teoBNP}
Assume that $F_V=0$ or $F_V $is a reduced connected divisor with normal crossings and contracted by $V$. Then the topological index $\tau(X)$ of $X$ is non-negative.
\end{theorem}
The condition on $F_V$ of being contracted by $V$ has been shown to be necessary, whereas the authors wonder whether the other assumptions can be dropped. The argument can be extended to reduced  non-connected $F_V$, giving a lower bound on the index that can be negative. The authors prove that the index is non-negative provided that $F_V$ has at most one connected component of arithmetic genus $0$.
As for the general cases, they conjecture the following.
\begin{conjecture}\label{conjecture GENERALIZED LAGRANGIAN}\cite[Conjecture 1]{BNP}
Let $X$ be a minimal generalized Lagrangian surface of general type, with $F_V$ contracted by $V$. Then $\tau (X)\geq 0$.
\end{conjecture}
\begin{conjecture}\label{conjecture LAGRANGIAN}\cite[Conjecture 2]{BNP}
Let $X$ be a Lagrangian surface. Then $\tau (X)\geq 0$.
\end{conjecture}
In this section we shall see that $\lpr$ surfaces provide counterexamples to both conjectures. As usual, let $(S, [\mathcal L])\in \mathcal W(1,2)$, and let $X$ be the $\lpr$ surface assoviated. By Theorems \ref{theorem TOPOLOGICAL INDEX} and \ref{theorem ALBANESE} we have that $X$ is a Lagrangian surface with negative topological index, hence it disproves Conjecture \ref{conjecture LAGRANGIAN}.

In particular, $X$ is a generalized Lagrangian surface with $V=H^0(X,\Omega^1_X)$ by Lemma \ref{lemma IRREGULARITY} and Theorem \ref{theorem IRREGULARITY}. We shall prove that for a general choice of $(S, [\mathcal L])\in \mathcal W(1,2)$, the divisor $F_V$ associated to $X$ is smooth, reduced and contracted by $V$, with $6$ connected rational components (see Theorem \ref{theorem EFFEVI}). Being $\tau (X)=-2$, this fact disproves Conjecture \ref{conjecture GENERALIZED LAGRANGIAN}, thus showing that the assumption of connectedness of $F_V$ in Theorem \ref{teoBNP} is necessary.

The estimate on the topological index in \cite{BNP} is obtained via a fine analysis of the image sheaf of the evaluation morphism
$V\otimes \oo_X\longrightarrow \Omega_X^1$.
The authors prove the following inequality (see the proof of \cite[Theorem 5.4]{BNP})
\begin{equation}\label{disuguaglianza EFFEVI}
\tau (X) \geq \frac{2}{3}K_X\cdot F_V+ \frac{1}{3}F_V^2.
\end{equation}
Thanks to the computation of $F_V$ in our examples, we shall see that this inequality is sharp. Indeed in our case the lower limit is negative, and it is achieved for a general choice of the abelian surface (see Remark \ref{remark SHARP} below).

\medskip
For $(S, [\mathcal L])\in \mathcal W(1,2)$, let us consider the corresponding $\lpr$ surface $X$ together with the degree two morphisms $\alpha_i\colon X\longrightarrow \overline{S}$, with $1\leq i\leq 3$. We set $V=H^0(X, \Omega_X^1)$ and we want to investigate the divisor $F_V$.

\begin{remark}
Let $a\colon X\longrightarrow \Alb(X)$ the Albanese map, which in this case is a morphism. Consider the dual of the differential of $a$:
$ da^*\colon H^0(X, \Omega_X^1)\longrightarrow \Omega^1_X, $ and its determinant
$$ \wedge^2 da^*\colon\bigwedge^2H^0(X, \Omega_X^1)\longrightarrow \omega_X .$$
The divisor $F_V$ is the divisorial part of the locus where $\wedge^2da^*$ has not maximal rank.
Note that even in this particular situation there could be curves not contracted by $a$ nevertheless belonging to $F_V$.
\end{remark}
Let us state the main result of the section.
Let $ \mathcal V, \mathcal U \subset \mathcal W(1,2)$ be the Zariski open subsets
$$
\mathcal V=\{ (S,[\mathcal L])\in   \mathcal W(1,2)\mid B_{\alpha_1} \mbox{ is smooth  } \},
$$
$$
\mathcal U=\{ (S,[\mathcal L])\in   \mathcal W(1,2)\mid B_{\alpha_1} \mbox{ is smooth at } B_{\alpha_1}\cap E_0\}.
$$
From Section \ref{section MONODROMY} we  know that these subsets are non-empty.

For $1\leq k\leq 3$, let us consider the rational curves $E_k\subset \overline{S}$ (cf. Notation \ref{notation ON barS}). We recall that any $E_k$ does not intersect the branch locus $\ba$, hence its pullback consists of two curves $E_k'$ and $E_k''$.

\begin{theorem}\label{theorem EFFEVI}
If $(S,[\mathcal L])\in \mathcal V$, then the divisor $F_V$ is reduced and consists of the $6$ smooth rational connected components $E_1',E_1'',\ldots, E_3''$, which are all $-3$-curves in $X$.

If $(S,[\mathcal L])\in \mathcal U$, then the divisor $F_V$ is reduced; its connected components are the $6$ smooth rational $-3$-curves $E_1',E_1'',\ldots, E_3''$, and possibly some smooth $-2$-curve.
\begin{proof}
Let us suppose that $\ba$ is smooth and  prove the first part of the statement. Under this assumption the Galois closure $W$ is smooth as well, and it coincides with $X$. Let $\beta \colon \overline S\longrightarrow S$ be the morphism constructed by blowing up the four base points of the linear system, $x_0,\ldots, x_3$, and then by blowing up the $6$ indeterminacy points $b_{jk}$ of the map $\gamma$, with $1\leq k\leq 3$ and $j=1,2$. As usual, we denote by $G_{jk}$ the exceptional divisor obtained by blowing up $b_{jk}$.

Let $u\in H^0(S, \omega_S)$ be the $2$-form on $S$. As $\alpha_i^*\beta^*u\in H^0(X, \omega_X)^{S_3}$ for any $1\leq i\leq 3$, the support of $F_V$ is contained in the intersection of the zero loci $Z(\alpha_i^*\beta^*u)$. By the formula for blow ups,
$$
K_{\overline S}= \beta^*K_S+ E_0+\sum_{k=1}^{3}E_k+2\sum_{k=1}^{3}(G_{1k}+G_{2k}).
$$
Hence $\beta^*u \in H^0(\overline{S}, \Omega^2_{\overline{S}})$ is such that $Z(\beta^*u)=E_0+\sum_{k=1}^{3}E_k+2\sum_{k=1}^{3}(G_{1k}+G_{2k})$.
The zero locus of $\alpha_i^*\beta^*u$ is $\alpha_i^*(Z(\beta^*u))+ R_{\alpha_i}$, where $R_{\alpha_i}$ denotes the ramification locus of $\alpha_i$. Looking at the definition of the $\alpha_i$'s, we have that
\begin{equation}\label{equation Ei}
\alpha_1^*(E_1+E_2+E_3)=\alpha_2^*(E_1+E_2+E_3)=\alpha_3^*(E_1+E_2+E_3).
\end{equation}
To see this fact, let us consider a point $p_1\in E_1$ and let $F$ be the fiber of $\overline{f}$ passing through $p_1$. Then $\overline{\gamma}^-1(\overline{\gamma}(p_1)=\{p_1,p_2,p_3\}$, where $p_2= E_2\cap F$ and $p_3= E_3\cap F$. Therefore
$$
\alpha_i^*(p_1+p_2+p_3)=\left\{(p_1,p_2),(p_1,p_3),(p_2,p_1),(p_2,p_3),(p_3,p_1),(p_3,p_2)\right\}
$$
for all $1\leq i\leq 3$ and (\ref{equation Ei}) follows. Analogously, by Lemma \ref{lemma TRIANGLE} and by the definition of the $\alpha_i$'s it is easy to check that the curves ${\alpha_i^*\left(E_0+2\sum_{k=1}^{3}(G_{1k}+G_{2k})\right)+R_{\alpha_i}}$ do not have common components. Thus we have that
$$
F_V\subset \bigcap_{i=1}^3 Z\left(\alpha_i^*\beta^*u\right)=\alpha_1^*(E_1+E_2+E_3)=E_1'+\ldots+E_3''.
$$
On the other hand, the pullbacks of the $E_k$'s are all contracted by the Albanese morphism $a$, because they are $-3$-curves of $X$. Thus $F_V=E_1'+\ldots+E_3''$ as claimed.

Suppose now that $W$ is singular. By Proposition \ref{proposition Balpha SINGULARITIES} it can only possess rational double points and - by the discussion in the previous section - these singularities can only lie either over $E_0\cap B_{\alpha_1}\cap B_{\alpha_2}\cap B_{\alpha_3}$, or over the images of the $80$  points of order three in $\overline S$.
In these cases, we can obtain the minimal desingularization of $W$ by resolving the singularities of the branch divisor $\ba$, and by taking a suitable double cover, the so-called canonical resolution (see \cite[Section III.7]{BHPV}). It is easy to check that the $-2$-curves arising in $X$ can appear as double components of $F_V$ only in the first case. The $-2$-curves corresponding to blowing up the order $3$ points are instead reduced in $F_V$. Hence, the statement is proved when $(S,[\mathcal L])\in \mathcal U$.
\end{proof}
\end{theorem}

\begin{remark}\label{remark SHARP}
Suppose that $(S,[\mathcal L])\in \mathcal V$, that is the Galois closure $W$ of $\overline{\gamma}$ is smooth and $W=X$. By the latter theorem we have that $F_V=\sum_{k=1}^3E_k'+E_k''$. Let us now consider the inequality (\ref{disuguaglianza EFFEVI}). By the adjunction formula, $(K_X+F_V)\cdot F_V=2 p_a(F_V)-2= -12$, whereas $F_V^{\cdot 2}=-18$. Hence the inequality (\ref{disuguaglianza EFFEVI}) gives $\tau (X)\geq -2$, which is actually an equality.
On the other hand, if $(S,[\mathcal L])\in \mathcal U$, the $\lpr$ surface is obtained from $W$ by blowing up the singular points and thus generating $-2$-curves. Let $n$ be the number of the $-2$-strings in the $\lpr$ surface $X$. Inequality (\ref{disuguaglianza EFFEVI}) becomes in this case $\tau (X)\geq -2-n$.
\end{remark}

\end{document}

%% file: FIBRAiperellittica.pstex_t
\begin{picture}(0,0)%
\includegraphics{FIBRAiperellittica.pstex}%
\end{picture}%
\setlength{\unitlength}{1579sp}%
\begingroup\makeatletter\ifx\SetFigFont\undefined%
\gdef\SetFigFont#1#2#3#4#5{%
  \reset@font\fontsize{#1}{#2pt}%
  \fontfamily{#3}\fontseries{#4}\fontshape{#5}%
  \selectfont}%
\fi\endgroup%
\begin{picture}(4682,4764)(7566,-3401)
\put(9335,1200){\makebox(0,0)[lb]{\smash{{\SetFigFont{11}{11.0}{\rmdefault}{\mddefault}{\updefault}{$G_{jk}$}%
}}}}
\put(10201,-1261){\makebox(0,0)[lb]{\smash{{\SetFigFont{8}{8.0}{\rmdefault}{\mddefault}{\updefault}{$2:1$}%
}}}}
\put(11026,-136){\makebox(0,0)[lb]{\smash{{\SetFigFont{8}{8.0}{\rmdefault}{\mddefault}{\updefault}{$1:1$}%
}}}}
\put(8101,-3337){\makebox(0,0)[lb]{\smash{{\SetFigFont{11}{11.0}{\rmdefault}{\mddefault}{\updefault}{$D$}%
}}}}
\put(12068,-3100){\makebox(0,0)[lb]{\smash{{\SetFigFont{11}{11.0}{\rmdefault}{\mddefault}{\updefault}{$\mathbb{P}^1$}%
}}}}
\put(7772,355){\makebox(0,0)[lb]{\smash{{\SetFigFont{11}{11.0}{\rmdefault}{\mddefault}{\updefault}{$p_0$}%
}}}}
\put(10557,630){\makebox(0,0)[lb]{\smash{{\SetFigFont{11}{11.0}{\rmdefault}{\mddefault}{\updefault}{$b_{jk}$}%
}}}}
\put(8300,-550){\makebox(0,0)[lb]{\smash{{\SetFigFont{11}{11.0}{\rmdefault}{\mddefault}{\updefault}{$\bar{b}_{jk}$}%
}}}}
\end{picture}%

%% file: FIBREbranch.pstex_t
\begin{picture}(0,0)%
\includegraphics{FIBREbranch.pstex}%
\end{picture}%
\setlength{\unitlength}{1184sp}%
\begingroup\makeatletter\ifx\SetFigFont\undefined%
\gdef\SetFigFont#1#2#3#4#5{%
  \reset@font\fontsize{#1}{#2pt}%
  \fontfamily{#3}\fontseries{#4}\fontshape{#5}%
  \selectfont}%
\fi\endgroup%
\begin{picture}(24719,9382)(-7221,-10218)
\put(-4567,-4103){\makebox(0,0)[lb]{\smash{{\SetFigFont{11}{11.0}{\rmdefault}{\mddefault}{\updefault}{$\textrm{(a)}$}%
}}}}
\put(5416,-2700){\makebox(0,0)[lb]{\smash{{\SetFigFont{11}{11.0}{\rmdefault}{\mddefault}{\updefault}{$p$}%
}}}}
\put(2145,-2131){\makebox(0,0)[lb]{\smash{{\SetFigFont{11}{11.0}{\rmdefault}{\mddefault}{\updefault}{$p_0$}%
}}}}
\put(14160,-3349){\makebox(0,0)[lb]{\smash{{\SetFigFont{11}{11.0}{\rmdefault}{\mddefault}{\updefault}{$p_0$}%
}}}}
\put(-6396,-1934){\makebox(0,0)[lb]{\smash{{\SetFigFont{11}{11.0}{\rmdefault}{\mddefault}{\updefault}{$p_0$}%
}}}}
\put(-1784,-1147){\makebox(0,0)[lb]{\smash{{\SetFigFont{11}{11.0}{\rmdefault}{\mddefault}{\updefault}{$q$}%
}}}}
\put(-4589,-1147){\makebox(0,0)[lb]{\smash{{\SetFigFont{11}{11.0}{\rmdefault}{\mddefault}{\updefault}{$p$}%
}}}}
\put(4762,-9668){\makebox(0,0)[lb]{\smash{{\SetFigFont{11}{11.0}{\rmdefault}{\mddefault}{\updefault}{$\textrm{(e)}$}%
}}}}
\put(-4567,-9668){\makebox(0,0)[lb]{\smash{{\SetFigFont{11}{11.0}{\rmdefault}{\mddefault}{\updefault}{$\textrm{(d)}$}%
}}}}
\put(13202,-4103){\makebox(0,0)[lb]{\smash{{\SetFigFont{11}{11.0}{\rmdefault}{\mddefault}{\updefault}{$\textrm{(c)}$}%
}}}}
\put(11550,-5930){\makebox(0,0)[lb]{\smash{{\SetFigFont{11}{11.0}{\rmdefault}{\mddefault}{\updefault}{$p_0$}%
}}}}
\put(15093,-7550){\makebox(0,0)[lb]{\smash{{\SetFigFont{11}{11.0}{\rmdefault}{\mddefault}{\updefault}{$q$}%
}}}}
\put(13202,-9668){\makebox(0,0)[lb]{\smash{{\SetFigFont{11}{11.0}{\rmdefault}{\mddefault}{\updefault}{$\textrm{(f)}$}%
}}}}
\put(8101,-8000){\makebox(0,0)[lb]{\smash{{\SetFigFont{11}{11.0}{\rmdefault}{\mddefault}{\updefault}{$q$}%
}}}}
\put(3601,-8000){\makebox(0,0)[lb]{\smash{{\SetFigFont{11}{11.0}{\rmdefault}{\mddefault}{\updefault}{$p_0$}%
}}}}
\put(-5851,-5668){\makebox(0,0)[lb]{\smash{{\SetFigFont{11}{11.0}{\rmdefault}{\mddefault}{\updefault}{$p_0$}%
}}}}
\put(-2515,-6210){\makebox(0,0)[lb]{\smash{{\SetFigFont{11}{11.0}{\rmdefault}{\mddefault}{\updefault}{$q$}%
}}}}
\put(4762,-4103){\makebox(0,0)[lb]{\smash{{\SetFigFont{11}{11.0}{\rmdefault}{\mddefault}{\updefault}{$\textrm{(b)}$}%
}}}}
\end{picture}%

%% file: FIBREtriangle.pstex_t
\begin{picture}(0,0)%
\includegraphics{FIBREtriangle.pstex}%
\end{picture}%
\setlength{\unitlength}{1184sp}%
\begingroup\makeatletter\ifx\SetFigFont\undefined%
\gdef\SetFigFont#1#2#3#4#5{%
  \reset@font\fontsize{#1}{#2pt}%
  \fontfamily{#3}\fontseries{#4}\fontshape{#5}%
  \selectfont}%
\fi\endgroup%
\begin{picture}(8655,10385)(3661,-3401)
\put(12301,1589){\makebox(0,0)[lb]{\smash{{\SetFigFont{11}{11.0}{\rmdefault}{\mddefault}{\updefault}{${\alpha_1}_{|D_2}$}%
}}}}
\put(8001,-3437){\makebox(0,0)[lb]{\smash{{\SetFigFont{11}{11.0}{\rmdefault}{\mddefault}{\updefault}{$D$}%
}}}}
\put(7650,300){\makebox(0,0)[lb]{\smash{{\SetFigFont{11}{11.0}{\rmdefault}{\mddefault}{\updefault}{$p_0$}%
}}}}
\put(9901,6089){\makebox(0,0)[lb]{\smash{{\SetFigFont{11}{11.0}{\rmdefault}{\mddefault}{\updefault}{$D_{2}$}%
}}}}
\put(7576,5289){\makebox(0,0)[lb]{\smash{{\SetFigFont{11}{11.0}{\rmdefault}{\mddefault}{\updefault}{$D_{3}$}%
}}}}
\put(10501,614){\makebox(0,0)[lb]{\smash{{\SetFigFont{11}{11.0}{\rmdefault}{\mddefault}{\updefault}{$G$}%
}}}}
\put(7876,3314){\makebox(0,0)[lb]{\smash{{\SetFigFont{11}{11.0}{\rmdefault}{\mddefault}{\updefault}{$D_{1}$}%
}}}}
\put(9976,1814){\makebox(0,0)[lb]{\smash{{\SetFigFont{9}{9.0}{\rmdefault}{\mddefault}{\updefault}{$2:1$}%
}}}}
\put(8326,1914){\makebox(0,0)[lb]{\smash{{\SetFigFont{11}{11.0}{\rmdefault}{\mddefault}{\updefault}{${\alpha_1}_{|D_1}$}%
}}}}
\put(3676,1589){\makebox(0,0)[lb]{\smash{{\SetFigFont{11}{11.0}{\rmdefault}{\mddefault}{\updefault}{${\alpha_1}_{|D_3}$}%
}}}}
\end{picture}%

%% file: lagrangian_14_giugno_2010_.bbl
\begin{thebibliography}{99}
\bibitem{wonderful} J. Amor{\'o}s, M. Burger,  K. Corlette, D.  Kotschick, D. Toledo,  \emph{Fundamental groups of compact {K}\"ahler manifolds}, {Mathematical Surveys and Monographs}, {Vol. 44}, AMS, {1996}.
\bibitem{ATV} M. Amram, M. Teicher, U. Vishne, The fundamental group of Galois cover of the surface $\mathbb T\times \mathbb T$, \emph{Internat. J. Algebra Comput.}, {\bf 18}(8) (2008), {1259--1282}.
\bibitem{BNP} M. A. Barja, J. C. Naranjo, G. P. Pirola, On the topological index of irregular surfaces, \emph{J. Algebraic Geom.}, {\bf 16}(4) (2007), 435--458.
\bibitem{Barth} W. Barth, Abelian surfaces with {$(1,2)$}-polarization, Algebraic geometry, Sendai, 1985,  Adv. Stud. Pure Math., {\bf 10}, 41--84, North-Holland, Amsterdam, 1987.
\bibitem{BHPV} W. Barth, K. Hulek, C. Peters, A. Van de Ven, \emph{Compact complex surfaces}, Second edition, Ergebnisse der Mathematik und ihrer Grenzgebiete {\bf 4}, Springer-Verlag, Berlin, 2004.
\bibitem{Beauville} A. Beauville, \emph{Complex algebraic surfaces},  {London Mathematical Society Lecture Note Series}, {Vol. 68}, {Cambridge University Press}, {1983}.
\bibitem{L-B} C. Birkenake, E. Lange, \emph{Complex abelian varieties}, Second edition, Grundlehren der Mathematischen Wissenschaften, {\bf 302}, Springer-Verlag, Berlin, 2004.
\bibitem{BT} F. Bogomolov, Y. Tschinkel, Lagrangian subvarieties of abelian fourfolds, \emph{Asian J. Math.}, {\bf 4}(1) (2000), 19--36.
\bibitem{Campana} F. Campana, Remarques sur les groupes de {K}\"ahler nilpotents, \emph{C. R. Acad. Sci. Paris S\'er. I Math.}, {\bf 317}(8), (1993), 777--780.
\bibitem{Catanese} F. Catanese, Moduli and classification of irregular K\"ahler manifolds (and algebraic varieties) with Albanese general type fibrations, \emph{Invent. Math.}, {\bf 104}(2) (1991), 263--289.
\bibitem{PC} A. Causin, G. P. Pirola,  Hermitian matrices and cohomology of {K}\"ahler varieties, \emph{Manuscripta Math.}, {\bf 121}(2) (2006),157--168.
\bibitem{CT} Z. Chen, S-L. Tan, Upper bounds in the Slope of a Genus 3 Fibration, \emph{Recent progress on some problems in several complex variables  and partial differential equations}, Contemp. Math., {\bf 400}, 65--87, AMS, 2006.
\bibitem{C-H} M. Cornalba, J. Harris, Divisor classes associated to families of stable varieties, with applications to the moduli space of curves, Ann. \emph{Sci. \'Ecole Norm. Sup.} (4) \textbf{21} (1988), 455--475.
\bibitem{DGMS}  P. Deligne, P. Griffiths, J. W. Morgan, D. Sullivan,  Real homotopy theory of {K}\"ahler manifolds, \emph{Invent. Math.}, {\bf 29}(3), (1975),245--274.
\bibitem{Faltings} G. Faltings, A new application of Diophantine approximations. \emph{A panorama of number theory or the view from Baker's garden (Z\"urich, 1999)}, 231--246, Cambridge Univ. Press, Cambridge, 2002.
\bibitem{Fujita} T. Fujita, On K\"aler fibre spaces over curves, \emph{J. Math. Soc. Japan}, {\bf 30}(4) (1978), 779--794.
\bibitem{FH} W. Fulton and J. Harris, \emph{Representation theory. A first course.}, Graduate Texts in Mathematics, {\bf 129}, Springer-Verlag, New York, 1991.
\bibitem{SGA} A. Grothendieck,  \emph{Rev\^etements \'etales et groupe fondamental. {F}asc. {I}},  {SGA}, {1960/61}, {Institut des Hautes \'Etudes Scientifiques},  {Paris}, {1963}.
\bibitem{Harris} J. Harris, Galois groups of enumerative problems, \emph{Duke Math. J.}, {\bf 46}(4) (1979), 685--724.
\bibitem{Hartshorne} R. Hartshorne, \emph{Algebraic geometry}, Graduate Texts in Mathematics, {\bf 52}, Springer-Verlag, New York-Heidelberg, 1977.
\bibitem{Hein} R. M. Hein, The geometry of the mixed Hodge structure on the fundamental group, \emph{Proc. Symposia in Pure Math.}, {\bf 46} (1987).
\bibitem{L-P} R. Lazarsfeld, M. Popa, BGG correspondence, cohomology of compact K\"ahler manifolds, and numerical inequalities, arXiv:0907.0651v1 [math.AG].
\bibitem{Liedtke} C. Liedtke, \emph{On fundamental groups of Galois closures of generic projections}, Bonner Mathematische Schriften, {\bf 367}, Dissertation, 2004.
\bibitem{Miyaoka} Y. Miyaoka Algebraic surfaces with positive indices, \emph{Classification of algebraic and analytic manifolds ({K}atata,1982)}, Progr. Math., {\bf 39}, Birkh\"auser Boston 1983, 281--301.

\bibitem{M-T} B. Moishezon, M. Teicher, Simply-connected algebraic surfaces of positive index, \emph{Invent. Math.},  {\bf 89}(3) (1987), {601--643}.
\bibitem{Morgan1} J. W. Morgan,  {\em The rational homotopy theory of smooth, complex projective varieties (following {P}. {D}eligne, {P}. {G}riffiths, {J}.  {M}organ, and {D}. {S}ullivan)} S\'eminaire {B}ourbaki, {V}ol. 1975/76, 28\`eme ann\'ee, {E}xp. {N}o. 475, 69--80. Lecture Notes in Math., Vol. 567. Springer, 1977.
\bibitem{Mumford}  D. Mumford, {\em Abelian varieties}, Tata Institute of Fundamental Research Studies in Mathematics, No. 5, Bombay, 1970.
\bibitem{Persson78} U. Persson, Double coverings and surfaces of general type, \emph{Algebraic geometry (Proc. Sympos., Univ. Troms\o, Troms\o, 1977)}, Lecture Notes in Mat., {\bf 687}, 168--195, Springer, Berlin, 1978.
\bibitem{Persson81} U. Persson, Chern invariants of surfaces of general type, \emph{Compositio Math.}, {\bf 43}(1) (1981), 3--58.
\bibitem{PIETRO} G. P. Pirola, Curves on generic {K}ummer varieties, \emph{Duke Math. J.}, {\bf 59}(3) (1989), 701--708.
\bibitem{Pirola} G. P. Pirola, On a conjecture of {X}iao, \emph{J. Reine Angew. Math.}, { \bf 431} (1992), 75--89.
\bibitem{Reid} M. Reid, \emph{Problems on pencils of small genus}, http://www.warwick.ac.uk/~masda/surf/
\bibitem{Siu} Y. T. Siu, Strong rigidity for {K}\"ahler manifolds and the construction of bounded holomorphic functions, \emph{Discrete groups in geometry and analysis ({N}ew {H}aven, 1984)}, {Progr. Math.}, 67, 124--151, Birkh\"auser Boston, 1987.
\bibitem{SVdV} A. J. Sommese, A. Van de Ven, Homotopy groups of pullbacks of varieties,\emph{ Nagoya Math. J.}, {\bf 102} (1986), 79--90.
\bibitem{TY} H. Tokunaga, H. Yoshihara, Degree of irrationality of abelian surfaces, \emph{J. Algebra}, {\bf 174}(3) (1995), 1111--1121.
\bibitem{Xiao} G. Xiao, Irregularity of surfaces with a linear pencil, \emph{Duke Math. J.}, {\bf 55}(3) (1987), 597--602.
\end{thebibliography}
